\documentclass[12pt,leqno]{article}
\usepackage{amssymb}
\usepackage[mathscr]{eucal}
\usepackage{amsmath,amssymb,latexsym,theorem,bbm}
\usepackage{amsmath,amssymb,latexsym,theorem,enumerate}
\usepackage{color,url}
\usepackage{appendix}

\newcommand{\DD}{\mathbb{D}}

\newcommand{\NN}{\mathbb{N}}

\newcommand{\RR}{\mathbb{R}}

\newcommand{\ZZ}{\mathbb{Z}}

\newcommand{\bX}{{\boldsymbol{X}}}

\newcommand{\bmu}{{\boldsymbol{\mu}}}

\newcommand{\bzero}{{\boldsymbol{0}}}

\newcommand{\cB}{{\mathcal B}}

\newcommand{\cD}{{\mathcal D}}

\newcommand{\cL}{{\mathcal L}}

\newcommand{\cN}{{\mathcal N}}

\newcommand{\cP}{{\mathcal P}}

\newcommand{\dd}{\mathrm{d}}
\newcommand{\ee}{\mathrm{e}}

\newcommand{\EE}{\operatorname{\mathbb{E}}}

\newcommand{\PP}{\operatorname{\mathbb{P}}}

\newcommand{\cov}{\operatorname{Cov}}

\renewcommand{\leq}{\leqslant}
\renewcommand{\geq}{\geqslant}

\newcommand{\stoch}{\stackrel{\PP}{\longrightarrow}}
\newcommand{\distr}{\stackrel{\cD}{\longrightarrow}}

\newcommand{\as}{\stackrel{{\mathrm{a.s.}}}{\longrightarrow}}

\newcommand{\bbone}{\mathbbm{1}}

\newcommand{\proofend}{\hfill\mbox{$\Box$}}

\numberwithin{equation}{section}

\theoremstyle{change} \theorembodyfont{\em}
\newtheorem{Lem}{Lemma.}[section]
\newtheorem{Thm}[Lem]{Theorem.}
\newtheorem{Pro}[Lem]{Proposition.}
\newtheorem{Cor}[Lem]{Corollary.}
\newtheorem{Def}[Lem]{Definition.}

\theorembodyfont{\rm}
\newtheorem{Rem}[Lem]{Remark.}

\begin{document}

\begin{center}
 {\bfseries\Large Limit theorems for Bajraktarevi\'c and Cauchy
  quotient means of independent identically distributed random variables}

\vspace*{3mm}

{\sc\large
  M\'aty\'as $\text{Barczy}^{*,\diamond}$,
  \ P\'al $\text{Burai}^{**}$ }

\end{center}

\vskip0.2cm

\noindent
 * MTA-SZTE Analysis and Stochastics Research Group,
   Bolyai Institute, University of Szeged,
   Aradi v\'ertan\'uk tere 1, H--6720 Szeged, Hungary.

\noindent
 ** Faculty of Informatics, University of Debrecen,
    Pf.~12, H--4010 Debrecen, Hungary.

\noindent e-mail: barczy@math.u-szeged.hu (M. Barczy),
                 burai.pal@inf.unideb.hu  (P. Burai).

\noindent $\diamond$ Corresponding author.

\vskip0.2cm


\renewcommand{\thefootnote}{}
\vspace*{0.2cm}
\footnote{M\'aty\'as Barczy is supported by the J\'anos Bolyai Research Scholarship
 of the Hungarian Academy of Sciences.}

\vspace*{-10mm}

\begin{abstract}
We derive strong  laws of large numbers and central limit theorems for Bajraktarevi\'c, Gini and exponential- (also called Beta-type) and
 logarithmic Cauchy quotient means of independent identically distributed (i.i.d.) random variables.
The exponential- and logarithmic Cauchy quotient means of a sequence of i.i.d.\ random variables
 behave asymptotically  normal with the usual square root scaling just like the geometric means of the given random variables.
Somewhat surprisingly, the multiplicative Cauchy quotient means of i.i.d.\ random variables behave asymptotically in a rather different way:
 in order to get a non-trivial normal limit distribution a time dependent centering is needed.
\end{abstract}

\section{Introduction}
\label{section_intro}

Studying properties  of various kinds of means (aggregation functions) is an old, popular and important topic due to the
  large number of applications in every branch of mathematics.
For a recent survey, see Beliakov et al.\ \cite{BelBusSan}.
This paper is devoted to studying the asymptotic behaviour of Bajraktarevi\'c means and Cauchy quotient means
 of independent identically distributed (i.i.d.) random variables.
Such an investigation for the arithmetic means of i.i.d. random variables goes back to Kolmogorov, and it is  at the heart of
 classical probability theory.
Recently, de Carvalho \cite[Theorem 1]{Car} (see also Theorem \ref{Thm_de_Carvalho}) has derived a central limit theorem for quasi arithmetic means,
 and he has also pointed out the fact that quasi arithmetic means have some applications in interest rate theory and
 unemployment duration analysis, see \cite[Examples 4 and 5]{Car}.

 We derive strong  laws of large numbers and central limit theorems for Bajraktarevi\'c, exponential Cauchy quotient and
 logarithmic Cauchy quotient means of i.i.d. random variables, see Theorems \ref{Thm_CLT_Baj_mean}, \ref{Thm_CLT_exp_Cauchy_quotient}
 and \ref{Thm_CLT_log_Cauchy_quotient}.
The multiplicative Cauchy quotient means of i.i.d.\ random variables behave asymptotically in a somewhat different way:
 in order to get a non-trivial normal limit distribution a time dependent centering is needed, see Theorem \ref{Thm_CLT_mult_Cauchy_quotient}.

We show another application of quasi arithmetic means to congressional apportionment in the USA's election
 motivated by Sullivan \cite{Sul1, Sul2}, and we also point out its possible extensions for Bajraktarevi\'c means and Cauchy quotient means,
 see Appendix \ref{Ex_Toredek_szavazatok}.%

Let \ $\NN$, \ $\ZZ_+$, \ $\RR$ \ and \ $\RR_+$ \ denote the sets of positive integers, non-negative integers, real numbers and non-negative real numbers.
Convergence almost surely, in probability and in distribution will be denoted by \ $\as$, \ $\stoch$ \ and \ $\distr$, \ respectively.
For any \ $d\in\NN$, \ $\cN_d(\bzero,\Sigma)$ \ denotes a \ $d$-dimensional normal distribution with mean vector \ $\bzero\in\RR^d$ \
 and covariance matrix \ $\Sigma\in\RR^{d\times d}$.
\ In  the case of \ $d=1$, \ instead of \ $\cN_1$ \ we simply write \ $\cN$.

\begin{Def}\label{Def_mean}
Let \ $I\subset \RR$ \ be an interval and \ $n\in\NN$.
\ A function \ $M:I^n \to \RR$ \ is called an \ $n$-variable mean in \ $I$ \ if
 \begin{align}\label{def_mean}
    \min(x_1,\ldots,x_n) \leq M(x_1,\ldots,x_n) \leq \max(x_1,\ldots,x_n),
     \qquad x_1,\ldots,x_n\in I.
 \end{align}
 An \ $n$-variable mean \ $M$ \ in \ $I$ \ is called strict if both inequalities in \eqref{def_mean} are sharp for all \ $x_1,\ldots,x_n\in I$ \
 satisfying \ $\min(x_1,\ldots,x_n)< \max(x_1,\ldots,x_n)$.
\end{Def}

If \ $n=1$, \ then the only \ $1$-variable mean \ $M$ \ in \ $I$ \ is \ $M(x)=x$, \ $x\in I$.

Kolmogorov and Nagumo provided an axiomatic construction for a sequence of functions \ $M_n:I^n \to \RR$, $n\in\NN$, \
 to define a  ''regular mean'' in \ $I$, \ where \ $I$ \ is a closed
 subinterval of \ $\RR$, \ see, e.g., Kolmogorov \cite{Kolmogoroff1930}, Nagumo \cite{Nag1, Nag2} and  Tikhomirov \cite[page 144]{Tik}.

\begin{Thm}[Kolmogorov  (1930) and Nagumo (1930)]\label{Thm_Kolmogorov}
Let \ $I$ \ be a closed and bounded subinterval of \ $\RR$, then the following two statements are equivalent:

\begin{enumerate}
\item[ (i)] There exists a sequence of functions $M_n:I^n\to\RR$, $n\in\NN$, such that
\begin{itemize}
  \item  $M_n$ \ is continuous and strictly monotone increasing in each variable for each \ $n\in\NN$,
  \item $M_n$ \ is symmetric for each \ $n\in\NN$ \  (i.e., \ $M_n(x_1,\ldots,x_n) = M_n(x_{\pi(1)},\ldots,x_{\pi(n)})$ \ for each
        \ $x_1,\ldots,x_n \in I$ \ and each permutation \ $(\pi(1),\ldots,\pi(n))$ \ of \ $(1,\ldots,n)$),
  \item $M_n(x_1,\ldots,x_n) = x$ \ whenever \ $x_1=\cdots =x_n = x\in I$, \ $n\in\NN$,
  \item $M_{n+m}(x_1,\ldots,x_n, y_1,\ldots,y_m) = M_{n+m}(\overline x_n, \ldots, \overline x_n, y_1,\ldots, y_m)$
        \ for each \ $n,m\in\NN$, \ $x_1,\ldots,x_n,y_1,\ldots,y_m\in I$, \ where \ $\overline x_n:=M_n(x_1,\ldots,x_n)$.
 \end{itemize}
\item[ (ii)] There exists a continuous and strictly monotone function \ $f:I\to\RR$ \ such that
 \[
   M_n(x_1,\ldots,x_n) = f^{-1} \left( \frac{1}{n} \sum_{i=1}^n f(x_i)\right),
    \qquad x_1,\ldots,x_n\in I,\;\; n\in\NN,
 \]
 where \ $f^{-1}$ \ denotes the inverse of \ $f$.
\end{enumerate}
\end{Thm}

\begin{Def}[Quasi arithmetic mean]
Let \ $n\in\NN$, \ let \ $I$ \ be a non-empty interval  of \ $\RR$, \ and let \ $f:I\to\RR$ \ be a continuous and strictly monotone increasing function.
The \ $n$-variable quasi arithmetic mean  of \ $x_1,\ldots, x_n\in I$ \ corresponding to \ $f$ \ is defined by
 \[
  M^f_n(x_1,\ldots,x_n):= f^{-1} \left( \frac{1}{n} \sum_{i=1}^n f(x_i)\right).
 \]
 The function \ $f$ \ is called a generator of \ $M^f_n$.
\end{Def}

\begin{Rem}
 (i).
The generator $f$ has an important role in the theory of quasi arithmetic means. It is not unique, but it is unique
 up to an affine transformation with nonzero factor (see, e.g.,  Hardy et al.\ \cite[Section 3.2, Theorem 83]{HLP1952}).
More precisely, two quasi arithmetic means  on \ $I$, \ generated by \ $f$ \ and \ $g$, \ are equal if and only if
 there  exist \ $a,b\in\RR$, \ $a\ne 0$ \ such that
 \begin{align}\label{help_quasi_generator}
  f(x)=ag(x)+b, \qquad x\in I.
 \end{align}
 As a consequence, the function \ $f$ \ in part (ii) of Theorem \ref{Thm_Kolmogorov} can be chosen to be strictly monotone increasing as well.

 (ii). A key idea in the theory of quasi arithmetic means is bisymmetry (for  the definition, see \eqref{E:bisymmetry} in Appendix \ref{App_non_quasi_arit_means}).
It is developed by  Acz\'el in \cite{Aczel1948}, who applied it for the characterization of \  $2$-variable quasi arithmetic means,
  and for the $n$-variable case, see M\"unnich et al.\ \cite{MunMakMok}.
The bisymmetry equation has importance also in the theory of quasisums and consistent aggregation in economical sciences (see, e.g.,
  Acz\'el and Maksa \cite{AczelMaksa1996}).

 (iii).  For more information about the story of quasi arithmetic means and their possible applications in   various areas,
 see the excellent survey of Muliere and Parmigiani \cite{MuliereParmigiani1993} and the references therein.
\proofend
\end{Rem}

For each \ $n\in\NN$, \ $M_n^f$ \ is a strict, symmetric \ $n$-variable mean in \ $I$ \ in the sense of Definition \ref{Def_mean}.
The arithmetic, geometric and harmonic mean is a quasi arithmetic mean corresponding to the function \ $f:\RR\to \RR$, $f(x):=x$, $x\in\RR$,
 \ $f:(0,\infty)\to \RR$, \ $f(x):=\ln(x)$, \ $x>0$, \ and \ \ $f:(0,\infty)\to \RR$, \ $f(x) =x^{-1}$, \ $x>0$, \ respectively.

Generalizing the notion of quasi arithmetic means, Bajraktarevi\'c \cite{Baj} introduced a new class of means
 (nowadays called Bajraktarevi\'c means) in the following way.

\begin{Def}[Bajraktarevi\'c mean]
Let \ $n\in\NN$, \ let \ $I$ \ be a non-empty interval  of \ $\RR$,
 let \ $f:I\to\RR$ \ be a continuous and strictly monotone
 function, and let \ $p:I\to (0,\infty)$ \ be a (weight) function.
The \ $n$-variable Bajraktarevi\'c mean  of \ $x_1,\ldots, x_n\in I$ \ corresponding to \ $f$ \ and \ $p$ \ is defined by
 \[
  B^{f,p}_n(x_1,\ldots,x_n):= f^{-1}\left( \frac{\sum_{i=1}^n p(x_i) f(x_i)}{\sum_{i=1}^n p(x_i)} \right).
 \]
\end{Def}

For each \ $n\in\NN$, \ $B^{f,p}_n$ \ is a strict, symmetric \ $n$-variable mean, see, e.g., Bajraktarevi\'c \cite{Baj} or P\'ales and Zakaria \cite{PalZak}.
Especially, by choosing \ $p(x)=1$, $x\in I$, \ we see that the Bajraktarevi\'c mean of \ $x_1,\ldots, x_n$ \ corresponding to \ $f$ \ and \ $p$ \
 coincides with the quasi arithmetic mean of \ $x_1,\ldots,x_n$ \ corresponding to \ $f$.

Next we recall the notion of Gini means, which are special Bajraktarevi\'c means.

\begin{Def}[Gini mean]\label{Def_Gini}
Let \ $r,s\in\RR$, \ $n\geq 2$, \ $n\in\NN$, \ and \ $x_1,\ldots,x_n>0$. \
The \ $n$-variable Gini mean of \ $x_1,\ldots, x_n$ \ corresponding to \ $r$ \ and \ $s$ \ is defined by
 \[
   G^{r,s}_n(x_1,\ldots,x_n)
    := \begin{cases}
        \left( \frac{ \sum_{i=1}^n x_i^r }{ \sum_{i=1}^n x_i^s } \right)^{\frac{1}{r-s}}  & \text{if \ $r\ne s$,}\\[2mm]
        \exp\left\{ \frac{ \sum_{i=1}^n x_i^s \ln(x_i) }{ \sum_{i=1}^n x_i^s }  \right\}
         = (\prod_{i=1}^n x_i^{x_i^s})^{\frac{1}{\sum_{i=1}^n x_i^s }}   & \text{if \ $r=s$.}
       \end{cases}
 \]
\end{Def}

Gini means are special Bajraktarevi\'c means, since, by choosing \ $I:=(0,\infty)$, \ $f:I\to\RR$,
 \[
  f(x):=\begin{cases}
          x^{\max(r,s)-\min(r,s)} & \text{if \ $r\ne s$,}\\
          \ln(x)  &  \text{if \ $r=s$,}
        \end{cases}
 \]
 and \ $p:I\to\RR$, \ $p(x):=x^{\min(r,s)}$, $x\in I$, \ the Bajraktarevi\'c mean of \ $x_1,\ldots, x_n\in I$ \ corresponding
  to $f$ and $p$ coincides with the Gini mean of $x_1,\ldots, x_n$ corresponding to $r$ and $s$.

Recently, Himmel and Matkowski \cite{HimMat2, HimMat3} have introduced and studied Cauchy quotient means.

\begin{Def}[Exponential Cauchy quotient mean, Beta-type mean]
Let \ $n\geq 2$, \ $n\in\NN$, \ and \ $x_1,\ldots,x_n>0$. \
The \ $n$-variable exponential Cauchy quotient mean  of \ $x_1,\ldots, x_n$ \ (also called \ $n$-variable Beta-type mean) is defined by
 \[
   \cB_n(x_1,\ldots,x_n):= \sqrt[n-1]{\frac{nx_1\cdots x_n}{x_1+\cdots+x_n}}.
 \]
\end{Def}

Note that \ $\cB_n$ \ is a strict,  symmetric \ $n$-variable mean in \ $(0,\infty)$ \ for each \ $n\geq 2$, \ $n\in\NN$, \ see Himmel and Matkowski \cite[Theorem 2]{HimMat2}.
In  the case of \ $n=2$, \ $\cB_n(x_1,x_2)$ \ coincides with the harmonic mean of \ $x_1$ and \ $x_2$, \ where \ $x_1,x_2>0$.

\begin{Def}[Logarithmic Cauchy quotient mean]
Let \ $n\geq 2$, \ $n\in\NN$, \ and \ $x_1,\ldots,x_n>1$. \
The \ $n$-variable logarithmic Cauchy quotient mean  of \ $x_1,\ldots, x_n$ \ is defined by
 \[
   \cL_n(x_1,\ldots,x_n):= \frac{\sum_{i=1}^n \sqrt[n-1]{\prod_{j=1,j\ne i}^n x_j}\ln(x_i)  }{\sum_{i=1}^n \ln(x_i)}.
 \]
\end{Def}

Note that \ $\cL_n$ \ is a strict,  symmetric \ $n$-variable mean in \ $(1,\infty)$ \ for each \ $n\geq 2$, \ $n\in\NN$, \ see Himmel and Matkowski \cite[Theorem 2]{HimMat3}.

\begin{Def}[Multiplicative (or power) Cauchy quotient mean]
Let \ $n\geq 2$, \ $n\in\NN$, \ and \ $x_1,\ldots,x_n>1$. \
The \ $n$-variable multiplicative (or power) Cauchy quotient mean  of \ $x_1,\ldots, x_n$ \ is defined by
 \[
   \cP_n(x_1,\ldots,x_n):= \left( \prod_{i=1}^n x_i^{\ln\left( \frac{\ln(x_1\cdots x_n)}{\ln(x_i)} \right)}\right)^{\frac{1}{n\ln(n)}}.
 \]
\end{Def}

Note that \ $\cP_n$ \ is a strict,  symmetric \ $n$-variable mean in \ $(1,\infty)$ \ for each \ $n\geq 2$, \ $n\in\NN$, \ see
 Appendix \ref{Appendix_mult_Cauchy_quotient_mean} or Himmel and Matkowski \cite[Theorem 2]{HimMat1}.
Since the reference \cite{HimMat1} refers to Himmel and Matkowski's slides of a talk given at a conference,
 where no proofs are available, and we have not found any other reference to the result in question,
 we decided to check that \ $\cP_n$ \ is indeed a strict \ $n$-variable mean in \ $(1,\infty)$ \ for each \ $n\geq 2$, \ $n\in\NN$,
 \ see Appendix \ref{Appendix_mult_Cauchy_quotient_mean}.

In the next remark we point out the fact that Bajraktarevi\'c means, and the considered Cauchy quotient means are not quasi
 arithmetic means in general.

\begin{Rem}
The class of Bajraktarevi\'c means strictly contains the class of quasi arithmetic means.
To see this, we check that not all the Gini means (as special Bajraktarevi\'c means) are quasi arithmetic means.
Gini means are trivially homogeneous, and a quasi arithmetic mean is homogeneous
 if and only if it is a H\"older mean  (also called power mean), i.e., it has the form
 \begin{align}\label{help_power_mean}
   \begin{cases}
     \left(\frac{1}{n}\sum_{i=1}^n x_i^p \right)^{\frac{1}{p}} & \text{if \ $p\ne 0$,}\\[1mm]
      \left( \prod_{i=1}^n x_i\right)^{\frac{1}{n}}   & \text{if \ $p=0$,}
   \end{cases}
   \qquad \forall\; x_1,\ldots,x_n>0,
 \end{align}
 with some \ $p\in\RR$, \ see, e.g., Hardy et al.\ \cite[ Section 3.3, Theorem 84]{HLP1952},
 and for some \ $n\geq 2$, \ $n\in\NN$, \ the class of \  $n$-variable Gini means strictly contains the class of \
  $n$-variable H\"older means (see, e.g.,  Bullen \cite[p. 248--251]{Bul}).

Himmel es Matkowski \cite[Remark 6]{HimMat2} showed that the exponential Cauchy quotient mean \ $\cB_n$ \ is a quasi arithmetic mean
 if and  only if \ $n=2$ \ (and in  the case of \ $n=2$, \ it is nothing else but the harmonic mean).
In Appendix \ref{App_non_quasi_arit_means}, we show that the logarithmic-, and multiplicative Cauchy quotient means \ $\cL_n$, \ $n\in\NN$, \
 and \ $\cP_n$, \ $n\in\NN$, \ are not quasi arithmetic means.
\proofend
\end{Rem}

 De Carvalho \cite[Theorem 1]{Car} derived a central limit theorem for quasi arithmetic means.
First, let us recall that if \ $f:I\to\RR$ \ is a continuous and strictly monotone increasing function,
 where \ $I$ \ is a non-empty subinterval of \ $\RR$, \ and \ $\xi$ \ is a random variable such that \ $\PP(\xi\in I)=1$ \ and
 \ $\EE(\vert f(\xi)\vert)<\infty$, \ then  Kolmogorov's expected value of \ $\xi$ \ corresponding to \ $f$ \ is defined by
 \[
   \EE_f(\xi):= f^{-1} (\EE (f(\xi)) ).
 \]
 Here \ $\EE (f(\xi))\in f(I)$, \ since \ $f(I)$ \ is an interval being a convex set.
If \ $I=(0,\infty)$ \ and \ $f(x) = x^p$, \ $x>0$, \ where \ $p>0$, \ then \ $\EE_f(\xi) = (\EE(\xi^p))^{\frac{1}{p}}$,
 \ which is nothing else, but the \ $L_p$-norm of \ $\xi$.
\ The usual expected value of \ $\xi$ \ corresponds to \ $f:\RR\to\RR$, \  $f(x):=ax+b$, \ $x\in\RR$, \ where \ $a,b\in\RR$, \ $a\ne 0$.
\ Recall also that \ $\DD^2(\xi):= \EE((\xi - \EE(\xi))^2)$ \ whenever \ $\EE(\vert\xi\vert)<\infty$.

\begin{Thm}[de Carvalho (2016)]\label{Thm_de_Carvalho}
Let \ $I$ \ be a non-empty  interval of \ $\RR$, \ and \ $f:I\to\RR$ \ be a continuous and strictly monotone increasing function.
Let \ $(\xi_n)_{n\in\NN}$ \ be a sequence of i.i.d.\ random variables such that
 \ $\PP(\xi_1\in I)=1$, \ $ \DD^2( f(\xi_1) )\in(0,\infty)$ \ and \ $f'(\EE_f(\xi_1))$ \  exists and  is non-zero.
Then
 \[
  M_n^f(\xi_1,\ldots,\xi_n) \as \EE_f(\xi_1) \qquad \text{as \ $n\to\infty$,}
 \]
 and
 \[
   \sqrt{n}\big(M_n^f(\xi_1,\ldots,\xi_n) - \EE_f(\xi_1)\big) \distr \cN\left(0, \frac{\DD^2(f(\xi_1))}{( f'(\EE_f(\xi_1)))^2} \right)
   \qquad \text{as \ $n\to\infty$.}
 \]
\end{Thm}

As a corollary of Theorem \ref{Thm_de_Carvalho}, de Carvalho \cite[Corollary 1]{Car} formulated central limit theorems for geometric and harmonic means.
We recall it for geometric means for our later purposes.

\begin{Cor}[de Carvalho (2016)]\label{Cor_de_Carvalho}
Let \ $(\eta_n)_{n\in\NN}$ \ be a sequence of i.i.d.\ random variables such that
 \ $\PP(\eta_1 > 0)=1$ \ and \ $ \DD^2( \ln(\eta_1) )\in(0,\infty)$.
\ Then
 \[
    \sqrt[n]{\eta_1\cdots\eta_n} \as \ee^{\EE(\ln(\eta_1))} \qquad \text{as \ $n\to\infty$,}
 \]
 and
 \begin{align}\label{help10_geom_mean}
   \sqrt{n} \big( \sqrt[n]{\eta_1\cdots\eta_n} - \ee^{\EE(\ln(\eta_1))} \big)
      \distr \cN\big(0, \DD^2(\ln(\eta_1))\ee^{2\EE(\ln(\eta_1))}  \big)
      \qquad \text{as \ $n\to\infty$.}
 \end{align}
\end{Cor}

Very recently, Mukhopadhyay et al.\ \cite[Lemma 3]{MukDasBasChaBha} have derived a central limit theorem for the power means (see \eqref{help_power_mean})
 of a sequence of independent random variables describing a mixture population consisting of two components: a major (dominating) and a minor (outlying) component.

 The paper is organized as follows.
Section \ref{section_results} contains our results, Section \ref{Section_proofs} is devoted to the proofs,
 and we close the paper with  four appendices, where we recall the Delta method (see Appendix \ref{App_delta_method}),
 we show that \ $\cP_n$ \ is a strict \ $n$-variable mean for each \ $n\geq 2$, \ $n\in\NN$ \ (see Appendix \ref{Appendix_mult_Cauchy_quotient_mean}),
 \ $\cL_n$ \ and  \ $\cP_n$ \ are not quasi arithmetic means for any \ $n\geq 2$, \ $n\in\NN$ \ (see Appendix \ref{App_non_quasi_arit_means}),
 and we give an application of quasi arithmetic and Bajraktarevi\'c means to congressional apportionment in the USA's election
 (see Appendix \ref{Ex_Toredek_szavazatok}).

\section{Results}
\label{section_results}

First, we present a strong law of large numbers and a central limit theorem for the Bajraktarevi\'c means
 of i.i.d.\ random variables.

\begin{Thm}\label{Thm_CLT_Baj_mean}
Let \ $I$ \ be a non-empty interval of \ $\RR$, \ let \ $f:I\to\RR$ \ be a continuous and strictly monotone
 function such that  the interval \ $f(I)$ \ is closed, and let \ $p:I\to (0,\infty)$ \ be a  measurable (weight) function.
Let \ $(\xi_n)_{n\in\NN}$ \ be a sequence of i.i.d.\ random variables such that \ $\PP(\xi_1\in I)=1$,
 \ $\EE((p(\xi_1))^2)<\infty$ \ and \ $\EE((p(\xi_1) f(\xi_1))^2)<\infty$.
\ If \ $f$ \ is differentiable at \ $f^{-1}\Big(\EE(p(\xi_1) f(\xi_1)) /  \EE(p(\xi_1)) \Big)$ \ with a non-zero derivative, then
 \[
    B^{f,p}_n(\xi_1,\ldots,\xi_n) \as f^{-1}\left( \frac{\EE(p(\xi_1)f(\xi_1))}{\EE(p(\xi_1))} \right)  \qquad \text{as \ $n\to\infty$,}
 \]
 and
 \begin{align}\label{help7_Baj_mean}
  \sqrt{n} \left(B^{f,p}_n(\xi_1,\ldots,\xi_n) - f^{-1}\left( \frac{\EE(p(\xi_1)f(\xi_1))}{\EE(p(\xi_1))} \right) \right) \distr \cN\big(0, \sigma^2_{f,p}\big)
    \qquad \text{as \ $n\to\infty$,}
 \end{align}
where
 \begin{align*}
  \sigma_{f,p}^2
    := \frac{(\EE(p(\xi_1)))^{-4}}
            {\left(f'\left( f^{-1}\left(\frac{\EE(p(\xi_1)f(\xi_1))}{\EE(p(\xi_1))}\right) \right)\right)^2}
         &\Big((\EE(p(\xi_1)))^2 \DD^2(p(\xi_1) f(\xi_1))  \\
         &\phantom{\Big(\;} - 2 \EE(p(\xi_1)) \EE(p(\xi_1)f(\xi_1)) \cov(p(\xi_1), p(\xi_1)f(\xi_1))\\
         &\phantom{\Big(\;}
             + (\EE(p(\xi_1) f(\xi_1) ))^2 \DD^2(p(\xi_1)) \Big).
 \end{align*}
\end{Thm}

Note that in Theorem  \ref{Thm_CLT_Baj_mean}, since \ $I$ \ is an interval and \ $f$ \ is continuous, we have \ $f(I)$ \
 \ is also an interval.
However, in general \ $f(I)$ \ is not closed, for example, if \ $I:=[0,\infty)$ \ and \ $f(x):=x/(x+1)$, $x\in I$, \ then \ $f(I)=[0,1)$. \
The assumption on the closedness of \ $f(I)$ \ in Theorem  \ref{Thm_CLT_Baj_mean} comes into play in proving
 a strong law of large numbers for \ $ B^{f,p}_n(\xi_1,\ldots,\xi_n)$ \ as \ $n\to\infty$.
\ Remark also that if \ $I=[a,b]$, \ where \ $a<b$, \ $a,b\in\RR$, \ and \ $f:I\to\RR$ \ is a continuous function,
 then \ $f(I)$ \ is closed, so in this special case the condition on the closedness of \ $f(I)$ \  in Theorem \ref{Thm_CLT_Baj_mean} is satisfied automatically.
One could easily specialize Theorem \ref{Thm_CLT_Baj_mean} for Gini means by choosing \ $f$ \ and \ $p$ \ as given
 after Definition \ref{Def_Gini}.

Next, we present a strong law of large numbers and a central limit theorem for the exponential Cauchy quotient means
 of i.i.d.\ random variables.

\begin{Thm}\label{Thm_CLT_exp_Cauchy_quotient}
Let \ $(\xi_n)_{n\in\NN}$ \ be a sequence of i.i.d.\ random variables such that
 \ $\PP(\xi_1 >0)=1$, \ $ \EE(\xi_1)<\infty$ \ and \ $\DD^2(\ln(\xi_1))\in(0,\infty)$.
\ Then
 \[
    \cB_n(\xi_1,\ldots,\xi_n) \as \ee^{\EE( \ln(\xi_1))}  \qquad \text{as \ $n\to\infty$,}
 \]
 and
 \begin{align}\label{help2_exp_Cauchy_quotient}
  \sqrt{n} \big(\cB_n(\xi_1,\ldots,\xi_n) - \ee^{\EE( \ln(\xi_1))}\big) \distr \cN\big(0, \DD^2(\ln(\xi_1)) \ee^{2\EE(\ln(\xi_1))}\big)
    \qquad \text{as \ $n\to\infty$.}
 \end{align}
\end{Thm}

\begin{Rem}
Concerning the moment conditions \ $\EE(\xi_1)<\infty$ \ and \ $\DD^2(\ln(\xi_1))\in(0,\infty)$ \ in Theorem \ref{Thm_CLT_exp_Cauchy_quotient},
 we note that they are not redundant in general.
Indeed, if \ $\xi_1:=\ee^{-\eta}$, \ where \ $\eta$ \ is a random variable such that \ $\PP(\eta\geq 0)=1$, \ $\EE(\eta)<\infty$ \ and
 \ $\EE(\eta^2)=\infty$, \ then \ $\PP(\xi_1>0)=1$, \ $\EE(\xi_1)\leq 1$, \ and \ $\EE((\ln(\xi_1))^2) = \EE(\eta^2)=\infty$.
\ Further, if \ $\xi_1:=\ee^{\eta}$, \ where \ $\eta$ \ is a random variable such that \ $\PP(\eta\geq 0)=1$, \ $\EE(\eta^2)<\infty$ \ and
 \ $\EE(\eta^3)=\infty$, \ then \ $\PP(\xi_1>0)=1$, \ $\EE((\ln(\xi_1))^2) = \EE(\eta^2)<\infty$, \ and
  \ $\EE(\xi_1) \geq \EE(\eta^3/3!) = \infty$ \ yielding that \ $\EE(\xi_1)=\infty$.
\proofend
\end{Rem}

Next, we present a strong law of large numbers and  central limit theorems for the logarithmic Cauchy quotient
 means of i.i.d.\ random variables.

\begin{Thm}\label{Thm_CLT_log_Cauchy_quotient}
Let \ $(\xi_n)_{n\in\NN}$ \ be a sequence of i.i.d.\ random variables such that
 \ $\PP(\xi_1 >1)=1$ \  and \ $\EE(\xi_1)<\infty$.
\ Then
 \[
    \cL_n(\xi_1,\ldots,\xi_n) \as \ee^{\EE( \ln(\xi_1))}
      \qquad \text{as \ $n\to\infty$,}
 \]
 and
 \begin{align}\label{help5_log_Cauchy_quotient}
  \sqrt{n} \big(\cL_n(\xi_1,\ldots,\xi_n) - \ee^{\EE( \ln(\xi_1))}\big) \distr \cN\big(0, \DD^2(\ln(\xi_1)) \ee^{2\EE(\ln(\xi_1))}\big)
    \qquad \text{as \ $n\to\infty$.}
 \end{align}
\end{Thm}

Note that the centralization \ $\ee^{\EE( \ln(\xi_1))}$ \ and the scaling \ $\sqrt{n}$ \ are the same
 in \eqref{help10_geom_mean}, \eqref{help2_exp_Cauchy_quotient} and in \eqref{help5_log_Cauchy_quotient},
 and the limit (normal) distributions coincide as well.
Roughly speaking, it means that the exponential- and logarithmic Cauchy quotient means of a sequence of i.i.d.\ random variables
 behave asymptotically just like the geometric means of the given random variables.

Next, we present a strong law of large numbers and a limit theorem for the multiplicative Cauchy quotient means of i.i.d.\ random variables.

\begin{Thm}\label{Thm_CLT_mult_Cauchy_quotient}
Let \ $(\xi_n)_{n\in\NN}$ \ be a sequence of i.i.d.\ random variables such that
 \ $\PP(\xi_1 >1)=1$ \  and \ $\DD^2(\ln(\xi_1))\in(0,\infty)$.
\begin{itemize}
 \item[(i)]
  Then
 \[
  \cP_n(\xi_1,\ldots,\xi_n) \as \ee^{\EE( \ln(\xi_1))}
      \qquad \text{as \ $n\to\infty$,}
 \]
 and
 \begin{align}\label{help5_mult_Cauchy_quotient}
  \ln(n) \big(\cP_n(\xi_1,\ldots,\xi_n) - \ee^{\EE( \ln(\xi_1))}\big) \stoch \ee^{\EE(\ln(\xi_1))}
           \Big( \ln(\EE(\ln(\xi_1))) \EE( \ln(\xi_1) ) - \EE(  \ln(\xi_1) \ln(\ln(\xi_1))) \Big)
 \end{align}
 as \ $n\to\infty$.

\item[(ii)]
  In addition, if \ $\DD^2(\ln(\xi_1) \ln(\ln(\xi_1)))\in(0,\infty)$, \ then
 \begin{align}\label{help5_mult_Cauchy_quotient_2}
 \begin{split}
 &\sqrt{n}\Bigg(
   \ln( \cP_n(\xi_1,\ldots,\xi_n) ) - \EE(\ln(\xi_1))
   - \frac{1}{\ln(n)}
      \Big( \ln(\EE(\ln(\xi_1))) \EE( \ln(\xi_1) ) - \EE(  \ln(\xi_1) \ln(\ln(\xi_1))) \Big)
   \Bigg)\\
 &\distr\cN(0, \DD^2(\ln(\xi_1)))
 \qquad \text{as \ $n\to\infty$.}
 \end{split}
 \end{align}
\end{itemize}
\end{Thm}

\begin{Rem}\label{Rem_P_n_momentumok}
Note that if \ $\PP(\xi_1 >1)=1$ \ and \ $\DD^2(\ln(\xi_1) \ln(\ln(\xi_1)))\in(0,\infty)$, \ then \ $\DD^2(\ln(\xi_1))\in(0,\infty)$.
\ Indeed,
 \begin{align*}
  \EE\left( (\ln(\xi_1))^2 \right)
    & = \EE\left( (\ln(\xi_1))^2 \bbone_{\{ \ln(\xi_1) \leq \ee\}}\right) +  \EE\left( (\ln(\xi_1))^2 \bbone_{\{ \ln(\xi_1) > \ee\} } \right) \\
    & \leq \ee^2 +  \EE\left( (\ln(\xi_1))^2 ( \ln(\ln(\xi_1)))^2 \bbone_{\{ \ln(\xi_1) > \ee\}}\right)
      \leq \ee^2+ \EE\left( (\ln(\xi_1))^2 ( \ln(\ln(\xi_1)))^2 \right)
     <\infty.
 \end{align*}
Next, we give an example  of a random variable \ $\xi_1$ \ such that
 \ $\PP(\xi_1 >1)=1$,  \ $\EE( (\ln(\xi_1))^2 )  < \infty$, \  and \ $\EE( (\ln(\xi_1))^2 (\ln(\ln(\xi_1)))^2  )  = \infty$,
 \ which shows that the condition \ $\DD^2(\ln(\xi_1) \ln(\ln(\xi_1)))\in(0,\infty)$ \ in part (ii) of Theorem \ref{Thm_CLT_mult_Cauchy_quotient}
 is indeed an additional one.
With the notation \ $\eta:=(\ln(\xi_1))^2$, \ it is enough to give an example  of a random variable \ $\eta$ \ such that
 \ $\PP(\eta\geq \ee)=1$, \ $\EE(\eta)<\infty$ \ and \ $\EE(\eta (\ln(\eta))^2)=\infty$.
\ Let \ $\eta$ \ be a random variable such that its density function takes the form
 \[
    f_\eta(x)=\begin{cases}
                 C\frac{1}{x^2(\ln(x))^2}  & \text{if \ $x\geq \ee$,}\\
                 0 & \text{if \ $x<\ee$,}
              \end{cases}
 \]
 where \ $\frac1C:=\int\limits_\ee^\infty \frac{1}{x^2(\ln(x))^2} \,\dd x$.
\ Note that \ $C\in(0,\infty)$, \ since with the substitution \ $x=\ee^y$,
 \[
  0<\int\limits_\ee^\infty \frac{1}{x^2(\ln(x))^2} \,\dd x
      = \int\limits_{1}^\infty \frac{1}{y^2\ee^y}\,\dd y
      \leq \int\limits_{1}^\infty \frac{1}{y^3}\,\dd y
      =\frac12.
 \]
Then \ $\PP(\eta\geq \ee)=1$, moreover,
 \begin{align*}
  \EE(\eta) = \int\limits_\ee^\infty xf_\eta(x)\,\dd x
            = C\int\limits_\ee^\infty \frac{1}{x(\ln(x))^2} \,\dd x
            = C\int\limits_1^\infty \frac{1}{y^2} \,\dd y
            =C<\infty,
 \end{align*}
 and
 \[
 \EE(\eta (\ln(\eta))^2)
     = \int\limits_\ee^\infty x(\ln(x))^2f_\eta(x)\,\dd x
     = C\int\limits_\ee^\infty \frac{1}{x} \,\dd x
     =\infty.
 \]
\proofend
\end{Rem}

Note that the limit distribution in \eqref{help5_mult_Cauchy_quotient} is not a normal distribution instead  of a deterministic constant,
 and the scaling factor is \ $\ln(n)$ \ instead of the usual \ $\sqrt{n}$.
So, somewhat surprisingly, the multiplicative Cauchy quotient means of i.i.d.\ random variables admit a different asymptotic behaviour than the exponential- and
 logarithmic Cauchy quotient means  of the random variables in question.

\section{Proofs}\label{Section_proofs}

\noindent{\bf Proof of Theorem \ref{Thm_CLT_Baj_mean}.}
By the strong law of large numbers,
 \begin{align}\label{help12}
   \frac{\frac{1}{n}\sum_{i=1}^n p(\xi_i) f(\xi_i)}{\frac{1}{n}\sum_{i=1}^n p(\xi_i)}
    \as \frac{\EE(p(\xi_1)f(\xi_1))}{\EE(p(\xi_1))}
    \qquad \text{as \ $n\to\infty$.}
 \end{align}
Since \ $I$ \ is an interval and \ $f$ \ is continuous, we have \ $f(I)$ \ is also an interval, yielding that
 \[
     \frac{\frac{1}{n}\sum_{i=1}^n p(\xi_i) f(\xi_i)}{\frac{1}{n}\sum_{i=1}^n p(\xi_i)}
        = \sum_{i=1}^n \frac{p(\xi_i)}{\sum_{j=1}^n p(\xi_j)} f(\xi_i) \in f(I),\qquad n\in\NN.
 \]
Using that \ $f(I)$ \ is assumed to be closed, by \eqref{help12}, we have
 \[
   \frac{\EE(p(\xi_1)f(\xi_1))}{\EE(p(\xi_1))} \in f(I),
 \]
 and hence, using that \ $f^{-1}$ \ is continuous,
 \begin{align*}
   B^{f,p}_n(\xi_1,\ldots,\xi_n)
    = f^{-1}\left( \frac{\frac{1}{n}\sum_{i=1}^n p(\xi_i) f(\xi_i)}{\frac{1}{n}\sum_{i=1}^n p(\xi_i)} \right)
    \as f^{-1}\left( \frac{\EE(p(\xi_1)f(\xi_1))}{\EE(p(\xi_1))} \right)  \qquad \text{as \ $n\to\infty$,}
 \end{align*}
 as desired.

By the multidimensional central limit theorem, we have
 \[
   \sqrt{n} \left(
            \begin{bmatrix}
              \frac{p(\xi_1)f(\xi_1)+\cdots+p(\xi_n)f(\xi_n)}{n} \\
              \frac{p(\xi_1)+\cdots+ p(\xi_n)}{n} \\
            \end{bmatrix}
            -
            \begin{bmatrix}
              \EE(p(\xi_1)f(\xi_1)) \\
              \EE(p(\xi_1)) \\
            \end{bmatrix}
            \right)
  \distr
    \cN_2   \left(  \begin{bmatrix}
              0 \\
              0 \\
            \end{bmatrix},
            \Sigma
             \right)
      \qquad \text{as \ $n\to\infty$,}
 \]
 where
 \[
  \Sigma:=\begin{bmatrix}
              \DD^2(p(\xi_1)f(\xi_1)) & \cov(p(\xi_1)f(\xi_1),p(\xi_1)) \\
              \cov(p(\xi_1)f(\xi_1),p(\xi_1)) & \DD^2(p(\xi_1)) \\
            \end{bmatrix}.
 \]
Using the Delta method with a measurable function \ $g:\RR^2\to\RR$ \ satisfying \ $g(x,y)=\frac{x}{y}$, \ $x,y>0$ \ (see, e.g., Theorem \ref{Thm_Delta_Method}),
 we have
 \begin{align*}
   \sqrt{n} \left(
               \frac{ \sum_{i=1}^n p(\xi_i)f(\xi_i) }{\sum_{i=1}^n p(\xi_i)}
               - \frac{\EE(p(\xi_1)f(\xi_1))}{\EE(p(\xi_1))}
            \right)
            \distr
            \cN(0,D\Sigma D^\top)
             \qquad \text{as \ $n\to\infty$,}
 \end{align*}
 where
 \[
   D:=g'( \EE(p(\xi_1)f(\xi_1)) , \EE(p(\xi_1)) )
       = \begin{bmatrix}
          \frac{1}{\EE(p(\xi_1))} & -\frac{\EE(p(\xi_1)f(\xi_1)) }{ (\EE(p(\xi_1)))^2}
         \end{bmatrix},
 \]
 and one can calculate
 \[
  D\Sigma D^\top
     =  \frac{\DD^2(p(\xi_1)f(\xi_1))}{(\EE(p(\xi_1)))^2}
        - 2 \cov( p(\xi_1)f(\xi_1), p(\xi_1)) \frac{\EE( p(\xi_1) f(\xi_1))}{(\EE(p(\xi_1)))^3}
        + \DD^2(p(\xi_1)) \frac{(\EE(p(\xi_1)f(\xi_1)))^2}{(\EE(p(\xi_1)))^4}.
 \]
Using again the Delta method with a measurable function \ $g:\RR^2\to\RR$ \ satisfying \ $g(x,y)=f^{-1}(x)$, \ $x\in f(I)$,
 we have
 \begin{align*}
   &\sqrt{n} \left(
               f^{-1}\left(\frac{ \sum_{i=1}^n p(\xi_i)f(\xi_i) }{\sum_{i=1}^n p(\xi_i)} \right)
               - f^{-1}\left( \frac{\EE(p(\xi_1)f(\xi_1))}{\EE(p(\xi_1))}\right)
            \right)\\
   &\qquad \distr \cN\left(0, \left( g'\left(\frac{\EE(p(\xi_1)f(\xi_1))}{\EE(p(\xi_1))} \right)\right)^2 D\Sigma D^\top\right)
      \qquad \text{as \ $n\to\infty$,}
 \end{align*}
  yielding the statement, since \ $g'(x) = 1/f'(f^{-1}(x))$, \ $x\in f(I)$.
\proofend

\noindent{\bf Proof of Theorem \ref{Thm_CLT_exp_Cauchy_quotient}.}
We have
 \begin{align}\label{help1}
    \cB_n(\xi_1,\ldots,\xi_n)
      = \sqrt[n-1]{\frac{n\xi_1\cdots\xi_n}{\xi_1+\cdots+\xi_n}}
      = \frac{ (\sqrt[n]{\xi_1\cdots\xi_n})^{\frac{n}{n-1}}}
             {\left( \frac{\xi_1+\cdots+\xi_n}{n} \right)^{\frac{1}{n-1}}},\qquad  n\geq 2, \;\; n\in\NN,
 \end{align}
 and hence the strong law of large numbers and Corollary \ref{Cor_de_Carvalho} yield that \ $\cB_n(\xi_1,\ldots,\xi_n)
 \as \frac{\ee^{\EE(\ln(\xi_1))}}{(\EE(\xi_1))^0} = \ee^{\EE(\ln(\xi_1))}$ \ as \ $n\to\infty$, \ as desired.

Further,
 \[
   \ln(\cB_n(\xi_1,\ldots,\xi_n))
     = \frac{1}{n-1} \sum_{i=1}^n \ln(\xi_i) - \frac{1}{n-1} \ln\left(\frac{\xi_1+\cdots+\xi_n}{n}\right) ,\qquad n\geq 2, \;\; n\in\NN,
 \]
 and the central limit theorem, Slutsky's lemma and \eqref{help1} yield that
 \begin{align*}
  &\sqrt{n}\big( \ln(\cB_n(\xi_1,\ldots,\xi_n)) - \EE(\ln(\xi_1)) \big) \\
  & = \sqrt{n} \left( \frac{n-1}{n} \ln(\cB_n(\xi_1,\ldots,\xi_n)) - \EE(\ln(\xi_1))
                      +\frac{1}{n} \ln(\cB_n(\xi_1,\ldots,\xi_n)) \right) \\
  & = \sqrt{n} \left(  \frac{1}{n} \sum_{i=1}^n \ln(\xi_i)  -  \EE(\ln(\xi_1))  \right)
       - \frac{1}{\sqrt{n}} \ln\left( \frac{\xi_1+\cdots+\xi_n}{n} \right)
       +\frac{1}{\sqrt{n}} \ln(\cB_n(\xi_1,\ldots,\xi_n)) \\
  & \distr \cN(0, \DD^2(\ln(\xi_1))) - 0\cdot \ln(\EE(\xi_1)) + 0\cdot \ln( \ee^{\EE(\ln(\xi_1))} )
     =  \cN(0, \DD^2(\ln(\xi_1)))
 \end{align*}
 as \ $n\to\infty$.
\ Using the  Delta method with the function \ $g:\RR\to\RR$, $g(x):=\ee^x$, \ $x\in\RR$ \ (see, e.g., Theorem \ref{Thm_Delta_Method}), we have
 \[
   \sqrt{n} ( \ee^{\ln(\cB_n(\xi_1,\ldots,\xi_n))} - \ee^{\EE( \ln(\xi_1))}) \distr \cN\big(0, \DD^2(\ln(\xi_1)) (\ee^{\EE(\ln(\xi_1))} )^2\big)
    \qquad \text{as \ $n\to\infty$,}
 \]
 yielding \eqref{help2_exp_Cauchy_quotient}.
\proofend

\noindent{\bf Proof of Theorem \ref{Thm_CLT_log_Cauchy_quotient}.}
 First note that, since \ $\PP(\ln(\xi_1)>0)=1$, \ we have
 \ $\xi_1 = \ee^{\ln(\xi_1)} \geq \frac{(\ln(\xi_1))^2}{2!}$ \ $\PP$-almost surely yielding that
 \ $\EE( (\ln(\xi_1))^2 ) \leq 2 \EE(\xi_1) <\infty$.

In the special case \ $\DD^2(\ln(\xi_1))=0$, \ we have \ $\PP(\xi_1=c)=1$ \ with some \ $c>1$, \ and
 \ $\cL_n(\xi_1,\ldots,\xi_n) = c$, $n\in\NN$, \ $\PP$-almost surely, yielding the assertion.
So in what follows, without loss of generality, we can assume that \ $\xi_1$ \ is non-degenerate, yielding that
 \ $\DD^2(\ln(\xi_1))\in(0,\infty)$.

For all \ $n\geq 2$, \ $n\in\NN$, \ we have
 \[
    \cL_n(\xi_1,\ldots,\xi_n)
      = \sqrt[n-1]{\prod_{j=1}^n \xi_j } \cdot
        \frac{\sum_{i=1}^n \xi_i^{-\frac{1}{n-1}} \ln(\xi_i) }
             {\sum_{i=1}^n \ln(\xi_i)},
 \]
 where, by Corollary \ref{Cor_de_Carvalho}, \ $\sqrt[n-1]{\prod_{j=1}^n \xi_j } =  \bigg(\sqrt[n]{\prod_{j=1}^n \xi_j }\bigg)^{\frac{n}{n-1}}
  \as \ee^{\EE(\ln(\xi_1))}$ \ as \ $n\to\infty$, \ and
   \begin{align}\label{help4}
     \frac{\sum_{i=1}^n \xi_i^{-\frac{1}{n-1}} \ln(\xi_i) }
          {\sum_{i=1}^n \ln(\xi_i)}
          \as 1 \qquad \text{as \ $n\to\infty$.}
  \end{align}
Indeed, since \ $\PP(\xi_1>1)=1$, \ we have \ $\PP(\xi_i^{-\frac{1}{n-1}}\ln(\xi_i)\leq \ln(\xi_i), i=1,\ldots,n)=1$, \ $n\geq 2$, \ $n\in\NN$, \ yielding
 \begin{align}\label{help8}
     \frac{\sum_{i=1}^n \xi_i^{-\frac{1}{n-1}} \ln(\xi_i) }
          {\sum_{i=1}^n \ln(\xi_i)}\leq 1,\qquad n\geq 2,
           \qquad \text{a.s.,}
 \end{align}
 and we also have
 \begin{align}\label{help9}
  \begin{split}
     \frac{\sum_{i=1}^n \xi_i^{-\frac{1}{n-1}} \ln(\xi_i) }
          {\sum_{i=1}^n \ln(\xi_i)}
         & \geq (\max(\xi_1,\ldots,\xi_n))^{-\frac{1}{n-1}}
           \geq (\xi_1+\cdots+\xi_n)^{-\frac{1}{n-1}}
           = \left(\!\frac{\xi_1+\cdots+\xi_n}{n}\!\right)^{-\frac{1}{n-1}} n^{-\frac{1}{n-1}} \\
         & \as (\EE(\xi_1))^0 \cdot 1 = 1
           \qquad \text{as \ $n\to\infty$.}
  \end{split}
 \end{align}
Consequently, by the squeeze theorem, we have \eqref{help4}, yielding that \ $\cL_n(\xi_1,\ldots,\xi_n)\as \ee^{\EE(\ln(\xi_1))}$ \
 as \ $n\to\infty$, \ as desired.

Further, for all \ $n\geq 2$,
 \begin{align*}
   &\sqrt{n} \big( \ln(\cL_n(\xi_1,\ldots,\xi_n)) - \EE(\ln(\xi_1)) \big)\\
   &\qquad  = \sqrt{n} \left( \frac{1}{n} \sum_{i=1}^n \ln(\xi_i) - \EE(\ln(\xi_1)) \right)
       + \frac{\sqrt{n}}{n(n-1)} \sum_{i=1}^n \ln(\xi_i)
       + \sqrt{n} \ln\left( \frac{\sum_{i=1}^n \xi_i^{-\frac{1}{n-1}} \ln(\xi_i) }{\sum_{i=1}^n \ln(\xi_i)} \right),
 \end{align*}
 where, by the central limit theorem,
 \begin{align}\label{help6}
   \sqrt{n} \left( \frac{1}{n} \sum_{i=1}^n \ln(\xi_i) - \EE(\ln(\xi_1)) \right) \distr \cN(0, \DD^2(\ln(\xi_1)))
    \qquad \text{as \ $n\to\infty$,}
 \end{align}
  and, by the strong law of large numbers,
  \[
    \frac{\sqrt{n}}{n-1}\cdot \frac{1}{n} \sum_{i=1}^n \ln(\xi_i) \as 0\cdot \EE(\ln(\xi_1)) = 0 \qquad \text{as \ $n\to\infty$,}
  \]
  and, by \eqref{help8}, \eqref{help9} and again the strong law of large numbers,
  \begin{align*}
    0&=\sqrt{n}\ln(1) \geq \sqrt{n} \ln\left( \frac{\sum_{i=1}^n \xi_i^{-\frac{1}{n-1}} \ln(\xi_i) }{\sum_{i=1}^n \ln(\xi_i)} \right)
                      \geq \sqrt{n}  \ln\left( \left(\frac{\xi_1+\cdots+\xi_n}{n}\right)^{-\frac{1}{n-1}} n^{-\frac{1}{n-1}}  \right)\\
                     & = -\frac{\sqrt{n}}{n-1} \ln\left( \frac{\xi_1+\cdots+\xi_n}{n}\right)
                         - \frac{\sqrt{n}}{n-1} \ln(n)
                      \as 0\cdot \ln(\EE(\xi_1)) - 0 = 0
         \qquad \text{as \ $n\to\infty$.}
   \end{align*}
Consequently, by Slutsky's lemma,
 \[
  \sqrt{n} \big( \ln(\cL_n(\xi_1,\ldots,\xi_n)) - \EE(\ln(\xi_1)) \big) \distr
    \cN(0, \DD^2(\ln(\xi_1))) \qquad \text{as \ $n\to\infty$,}
 \]
 and, an application of the Delta method (see, e.g., Theorem \ref{Thm_Delta_Method}) with the function \ $g:\RR\to\RR$, \ $g(x):=\ee^x$, \ $x\in\RR$, \
 yields \eqref{help5_log_Cauchy_quotient}.
\proofend

\noindent{\bf Proof of Theorem \ref{Thm_CLT_mult_Cauchy_quotient}.}
 First note that \ $\EE( \vert \ln(\xi_1) \ln(\ln(\xi_1))\vert)<\infty$.
\ Indeed, \ $\PP(\ln(\xi_1) > 0)=1$, \ and using that \ $1-\frac{1}{x} \leq \ln(x) \leq x-1$, \ $x>0$,
 \ we have
 \[
   \vert x\ln(x)\vert \leq \max(x^2+x,x+1)\leq x^2+2x+1=(x+1)^2,\qquad x>0,
 \]
 yielding that \ $\EE(\vert \ln(\xi_1) \ln(\ln(\xi_1))\vert )\leq \EE((\ln(\xi_1)+1)^2) <\infty$.

 (i). For all \ $n\geq 2$, \ $n\in\NN$, \ we have
 \begin{align}\label{help13}
 \begin{split}
  \ln(\cP_n(\xi_1&,\ldots,\xi_n))
    = \frac{1}{n\ln(n)} \sum_{i=1}^n \ln \left( \xi_i^{\ln\left( \frac{\ln(\xi_1\cdots \xi_n)}{\ln(\xi_i)} \right)}\right)\\
   & = \frac{1}{n\ln(n)} \sum_{i=1}^n \ln\left( \frac{\ln(\xi_1\cdots \xi_n)}{\ln(\xi_i)} \right) \ln(\xi_i) \\
   & = \frac{1}{n\ln(n)} \left( \ln(\ln(\xi_1\cdots \xi_n)) \sum_{i=1}^n \ln(\xi_i)
                                - \sum_{i=1}^n \ln(\xi_i)\ln(\ln(\xi_i))
                        \right)\\
   & = \frac{1}{\ln(n)} \ln\left(  \sum_{i=1}^n \ln(\xi_i) \right)
       \cdot \frac{1}{n}  \sum_{i=1}^n \ln(\xi_i)
      - \frac{1}{\ln(n)} \cdot \frac{1}{n} \sum_{i=1}^n \ln(\xi_i) \ln(\ln(\xi_i))\\
   & = \left( 1+ \frac{\ln\left(  \frac{1}{n}\sum_{i=1}^n \ln(\xi_i) \right)}{\ln(n)} \right)
       \frac{1}{n}  \sum_{i=1}^n \ln(\xi_i)
       - \frac{1}{\ln(n)} \cdot \frac{1}{n} \sum_{i=1}^n \ln(\xi_i) \ln(\ln(\xi_i)).
 \end{split}
 \end{align}
Hence, by the strong law of large numbers,
 \begin{align*}
  \ln(\cP_n(\xi_1,\ldots,\xi_n))
   & \as \left( 1+ 0\cdot\ln(\EE(\ln(\xi_1))) \right) \EE(\ln(\xi_1))
         - 0\cdot \EE(\ln(\xi_1) \ln(\ln(\xi_1)))
    = \EE(\ln(\xi_1))
 \end{align*}
 as \ $n\to\infty$, \ yielding \ $\cP_n(\xi_1,\ldots,\xi_n) \as \ee^{\EE(\ln(\xi_1))}$ \ as \ $n\to\infty$, \ as desired.

Further, by  the strong law of large numbers, \eqref{help6} and Slutsky's lemma,
 \begin{align*}
   &\ln(n) \left( \ln(\cP_n(\xi_1,\ldots,\xi_n)) - \EE(\ln(\xi_1)) \right) \\
   &  = \frac{\ln(n)}{\sqrt{n}} \cdot \sqrt{n} \left( \frac{1}{n} \sum_{i=1}^n \ln(\xi_i)  - \EE(\ln(\xi_1)) \right)
        + \ln \left(  \frac{1}{n}\sum_{i=1}^n \ln(\xi_i) \right) \frac{1}{n} \sum_{i=1}^n \ln(\xi_i) \\
   &\phantom{=\;}   - \frac{1}{n} \sum_{i=1}^n \ln(\xi_i) \ln(\ln(\xi_i)) \\
   &\distr  \ln(\EE(\ln(\xi_1))) \EE(\ln(\xi_1))  - \EE(\ln(\xi_1) \ln(\ln(\xi_1)))
   \qquad \text{as \ $n\to\infty$.}
 \end{align*}
Since the limit \ $\ln(\EE(\xi_1)) \EE(\ln(\xi_1))  - \EE(\ln(\xi_1) \ln(\ln(\xi_1)))$ \ is a constant, we also have
 \[
   \ln(n) \left( \ln(\cP_n(\xi_1,\ldots,\xi_n)) - \EE(\ln(\xi_1)) \right)
    \stoch   \ln(\EE(\ln(\xi_1))) \EE(\ln(\xi_1))  - \EE(\ln(\xi_1) \ln(\ln(\xi_1)))
 \]
 as \ $n\to\infty$.
\ Finally, an application of the Delta method (see, e.g., Theorem \ref{Thm_Delta_Method}) with the function
 \ $g:\RR\to\RR$, \ $g(x):=\ee^x$, \ $x\in\RR$, \ yields that
 \[
  \ln(n) (\cP_n(\xi_1,\ldots,\xi_n) - \ee^{\EE( \ln(\xi_1))}) \distr \ee^{\EE(\ln(\xi_1))}
           \big( \ln(\EE(\xi_1)) \EE( \ln(\xi_1) ) - \EE(  \ln(\xi_1) \ln(\ln(\xi_1))) \big)
 \]
 as \ $n\to\infty$.
Using that the limit \ $\ee^{\EE(\ln(\xi_1))}\big(  \ln(\EE(\ln(\xi_1))) \EE( \ln(\xi_1) ) - \EE(  \ln(\xi_1) \ln(\ln(\xi_1))) \big)$ \ is a constant,
 we have \eqref{help5_mult_Cauchy_quotient}, as desired.

(ii). First recall that \ $\DD^2(\ln(\xi_1))\in(0,\infty)$, \ see Remark \ref{Rem_P_n_momentumok}.
Using \eqref{help13}, for each \ $n\in\NN$, \ we have
 \begin{align*}
 &\sqrt{n}\Bigg(
   \ln( \cP_n(\xi_1,\ldots,\xi_n) ) - \EE(\ln(\xi_1))
   - \frac{1}{\ln(n)}
      \Big( \ln(\EE(\ln(\xi_1))) \EE( \ln(\xi_1) ) - \EE(  \ln(\xi_1) \ln(\ln(\xi_1))) \Big)
   \Bigg) \\
 &\qquad  = A^{(1)}_n +  A^{(2)}_n + A^{(3)}_n + A^{(4)}_n,
 \end{align*}
 where
 \begin{align*}
  A^{(1)}_n&:=\sqrt{n} \left( \frac{1}{n}  \sum_{i=1}^n \ln(\xi_i) - \EE( \ln(\xi_1) ) \right) ,\\
  A^{(2)}_n&:=-\frac{1}{\ln(n)} \sqrt{n} \left( \frac{1}{n} \sum_{i=1}^n \ln(\xi_i) \ln(\ln(\xi_i)) - \EE( \ln(\xi_1) \ln(\ln(\xi_1)) )  \right), \\
  A^{(3)}_n&:=\frac{1}{\ln(n)} \sqrt{n} \left(  \ln\left(  \frac{1}{n}\sum_{i=1}^n \ln(\xi_i) \right)  -  \ln(\EE(\ln(\xi_1))) \right)
                               \frac{1}{n} \sum_{i=1}^n \ln(\xi_i),\\
 A^{(4)}_n&:= \frac{1}{\ln(n)} \ln(\EE(\ln(\xi_1))) \sqrt{n} \left( \frac{1}{n}  \sum_{i=1}^n \ln(\xi_i) - \EE( \ln(\xi_1) ) \right).
 \end{align*}
To prove \eqref{help5_mult_Cauchy_quotient_2}, by Slutsky's lemma, it is enough to check that
 \ $A^{(1)}_n \distr \cN(0, \DD^2(\ln(\xi_1)))$ \ as \ $n\to\infty$, \ and \ $A^{(i)}_n\stoch 0$ \ as \ $n\to\infty$, \ $i=2,3,4$.
\ By the central limit theorem,
 \begin{align*}
   A^{(1)}_n \distr \cN(0, \DD^2(\ln(\xi_1))) \qquad \text{as \ $n\to\infty$,}
 \end{align*}
 and
 \begin{align*}
   \sqrt{n} \left( \frac{1}{n} \sum_{i=1}^n \ln(\xi_i) \ln(\ln(\xi_i)) - \EE( \ln(\xi_1) \ln(\ln(\xi_1)) )  \right)
          \distr \cN\big(0,\DD^2(\ln(\xi_1) \ln(\ln(\xi_1)))\big)
 \end{align*}
 as \ $n\to\infty$.
Hence, using Slutsky's lemma, we have \ $A^{(2)}_n\stoch 0$ \ as \ $n\to\infty$, \ and, using also that $A^{(4)}_n=A^{(1)}_n\frac{1}{\ln(n)} \ln(\EE(\ln(\xi_1)))$,
 we have \ $A^{(4)}_n\stoch 0$ \ as \ $n\to\infty$.

It remains to check that \ $A^{(3)}_n\stoch 0$ \ as \ $n\to\infty$.
 \ An application of the Delta method (see, e.g., Theorem \ref{Thm_Delta_Method}) with a measurable function \ $g:\RR\to\RR$ \ satisfying
\ $g(x)=\ln(x)$, \ $x>0$, \ yields that
 \[
    \sqrt{n} \left( \ln\left( \frac{1}{n}  \sum_{i=1}^n \ln(\xi_i) \right) - \ln(\EE( \ln(\xi_1))) \right)
      \distr \cN\left(0, \frac{\DD^2(\ln(\xi_1))}{( \EE(\ln(\xi_1)))^2}\right)
      \qquad \text{as \ $n\to\infty$.}
 \]
By the strong law of large numbers, we have \ $\frac{1}{n} \sum_{i=1}^n \ln(\xi_i) \as \EE(\ln(\xi_1))$ \ as \ $n\to\infty$.
Consequently, by Slutsky's lemma, we have \ $A^{(3)}_n\stoch 0$ \ as \ $n\to\infty$, \ as desired.
\proofend

\appendix

\vspace*{5mm}

\noindent{\bf\Large Appendices}

\section{Delta method}\label{App_delta_method}

We recall the  Delta method which we use for proving limit theorems, especially asymptotic normality, see, e.g.,
 Lehmann and Romano \cite[Theorem 11.2.14]{LehRom}.

\begin{Thm}\label{Thm_Delta_Method}
Let \ $\bX_n$, $n\in\NN$, \ and \ $\bX$ \ be \ $d$-dimensional random variables, where \ $d\in\NN$.
\ Assume that \ $\tau_n(\bX_n -\bmu) \distr \bX$ \  as \ $n\to\infty$ \ with some \ $\bmu\in\RR^d$ \ and \ $\tau_n\in\RR$, $n\in\NN$, \ satisfying
 \ $\tau_n\to\infty$ \ as \ $n\to\infty$.
\begin{enumerate}
  \item[(i)] Let \ $g:\RR^d\to\RR$ \ be a measurable function which is differentiable at \ $\bmu$.
        \ Then
        \[
          \tau_n (g(\bX_n) - g(\bmu)) \distr g'(\bmu)\bX \qquad \text{as \ $n\to\infty$,}
        \]
        where the $1\times d$ \ matrix \ $g'(\bmu)$  \ denotes the derivative of \ $g$ \ at \ $\bmu$.
        \ In particular, if \ $\bX$ \ is a \ $d$-dimensional normally distributed random variable with mean vector
        \ $\bzero\in\RR^d$ \ and covariance matrix \ $\Sigma\in\RR^{d\times d}$, \ then
        \[
         \tau_n (g(\bX_n) - g(\bmu)) \distr \cN(\bzero, g'(\bmu)\Sigma g'(\bmu)^\top)
          \qquad \text{as \ $n\to\infty$.}
        \]
  \item[(ii)]  More generally, let \ $g = (g_1,\ldots,g_q)^\top:\RR^d\to\RR^q$ \ be a measurable function
        which is differentiable at \ $\bmu$, \ where \ $d,q\in\NN$.
        \ Then
        \[
          \tau_n(g(\bX_n) - g(\bmu))
            = \tau_n ( g_1(\bX_n) - g_1(\bmu), g_2(\bX_n) - g_2(\bmu), \ldots, g_q(\bX_n) - g_q(\bmu) )^\top
              \distr g'(\bmu)\bX
        \]
        as \ $n\to\infty$, \
        where the $q\times d$ matrix \ $g'(\bmu)$  \ denotes the derivative of \ $g$ \ at \ $\bmu$.
        \ In particular, if \ $\bX$ \ is a \ $d$-dimensional normally distributed random variable with mean vector
        \ $\bzero\in\RR^d$ \ and covariance matrix \ $\Sigma\in\RR^{d\times d}$, \ then
        \[
         \tau_n (g(\bX_n) - g(\bmu)) \distr \cN_q(\bzero, g'(\bmu)\Sigma g'(\bmu)^\top)
          \qquad \text{as \ $n\to\infty$.}
        \]
\end{enumerate}
\end{Thm}

\section{$\cP_n$ \ is a strict \ $n$-variable mean} \label{Appendix_mult_Cauchy_quotient_mean}

\begin{Pro}
For each \ $n\geq 2$, \ $n\in\NN$, \ the multiplicative Cauchy quotient mean \ $\cP_n$ is a strict \ $n$-variable mean in \ $(1,\infty)$.
\end{Pro}

\noindent{\bf Proof.}
Let \ $x_1,\ldots, x_n>1$ \ be fixed such that \ $x_1\leq x_2\leq \cdots\leq x_n$.
With the notation \ $y_i:=\ln(x_i)$, \ $i=1,\ldots,n$, \ we have  that \ $ \min(x_1,\ldots,x_n) \leq \cP_n(x_1,\ldots,x_n) \leq \max(x_1,\ldots,x_n)$
 \ is equivalent to
 \[
   x_1^{n\ln(n)} \leq x_1^{\ln\left( \frac{\ln(x_1\cdots x_n)}{\ln(x_1)} \right)} \cdots x_n^{\ln\left( \frac{\ln(x_1\cdots x_n)}{\ln(x_n)} \right)}
                 \leq x_n^{n\ln(n)},
 \]
 which is equivalent to
 \[
   \ee^{n\ln(n) y_1} \leq \ee^{ y_1 \ln\left( \frac{y_1+\cdots + y_n}{y_1} \right) } \cdots \ee^{ y_n \ln\left( \frac{y_1+\cdots + y_n}{y_n} \right) }
                     \leq  \ee^{n\ln(n) y_n} ,
 \]
 or equivalently
 \begin{align}\label{help11}
   n\ln(n) y_1 \leq \sum_{i=1}^n y_i \ln\left( \frac{y_1+\cdots + y_n}{y_i} \right) \leq n\ln(n) y_n.
 \end{align}
Since \ $\ln\left( \frac{y_1+\cdots + y_n}{y_i} \right) = \ln\left( 1+ \frac{y_1+\cdots + y_{i-1} + y_{i+1}+ \cdots + y_n}{y_i} \right)\geq 0$, \
 and \ $0<y_1\leq y_2\leq \cdots \leq y_n$, \ we have
 \begin{align*}
   \sum_{i=1}^n y_i \ln\left( \frac{y_1+\cdots + y_n}{y_i} \right)
     \geq y_1 \sum_{i=1}^n \ln\left( \frac{y_1+\cdots + y_n}{y_i} \right),
 \end{align*}
 so for \ $n\ln(n) y_1 \leq \sum_{i=1}^n y_i \ln\left( \frac{y_1+\cdots + y_n}{y_i} \right)$ \ it is enough to check that
 \[
   n\ln(n) y_1 \leq y_1 \sum_{i=1}^n \ln\left( \frac{y_1+\cdots + y_n}{y_i} \right),
 \]
 or equivalently,
 \[
     \ln(n^n) \leq \ln\left( \frac{(y_1+\cdots + y_n)^n}{y_1\cdots y_n} \right).
 \]
By algebraic calculations, it is equivalent to \ $\sqrt[n]{y_1\cdots y_n}\leq (y_1+\cdots +y_n)/n$, \ which is nothing else but the
 well-known inequality between the arithmetic and geometric means, yielding that the first inequality in \eqref{help11} holds.

Now we turn to prove the second inequality in \eqref{help11}.
With the notation \ $z_i:=\frac{y_i}{y_n}$, \ $i=1,\ldots,n-1$, \  after dividing by \ $y_n$, \  we get that the
 second inequality in \eqref{help11} is equivalent to
 \[
   f(z_1,\ldots,z_{n-1}):=
    \sum_{i=1}^{n-1} z_i \ln\left( \frac{1+z_1+ \cdots + z_{n-1} }{z_i} \right)
      + \ln(z_1+\cdots + z_{n-1} + 1)
    \leq n\ln(n)
 \]
 for each \ $z_i\in(0,1]$, \ $i=1,\ldots, n-1$.
We check that the function \ $f: (0,1]^{n-1}\to\RR$ \ is strictly monotone increasing in each of its variables.
Due to the fact that \ $f$ \ is symmetric, it is enough to check it for the its first variable \ $z_1$.
\ One can calculate that
 \[
  \frac{\partial f}{\partial z_1} (z_1,\ldots,z_{n-1})
     = \ln\left( 1+ \frac{1+z_2+ \cdots + z_{n-1} }{z_1} \right) > 0,
     \qquad z_i\in(0,1], \;\; i=1,\ldots,n-1,
 \]
 yielding that \ $f$ \ is strictly monotone increasing in \ $z_1$.
\ Further, \ $f$ \ can be extended continuously onto \ $[0,1]^{n-1}$,
 \ since for any \ $a\in\RR_+$, \ by L'Hospital's rule,
 \[
   \lim_{x\downarrow 0} x \ln\left(1 + \frac{1+a}{x}\right) = \lim_{x\downarrow 0} \frac{\frac{1}{1+(1+a)/x} \frac{1+a}{x^2}}{1/x^2}=0.
 \]
Consequently, the function \ $f$ \ takes its maximum at \ $(1,\ldots,1)^\top\in\RR^{n-1}$, \ and
 \ $f(1,\ldots,1) = n\ln(n)$, \ yielding the second inequality in \eqref{help11}.

Finally, we present another proof of the second inequality in \eqref{help11}.
With the notation
 \[
  p_i:=\frac{y_i}{y_1+\cdots+y_n}\in(0,1),\qquad i = 1,\dots,n,
 \]
 the second inequality in \eqref{help11} takes the form
 \[
 -p_1\ln(p_1)-\cdots - p_n\ln(p_n) \leq n \ln(n) p_n.
 \]
Recall that \ $-p_1\log_2(p_1)-\cdots - p_n\log_2(p_n)$ \ is the entropy of the probability distribution \ $\{p_1,\ldots,p_n\}$,
 \ and it is well-known that the entropy of a probability distribution  concentrated  at \ $n$ \ points at most
 is less than or equal to \ $\log_2(n)$, \ yielding that
 \[
 -p_1\ln(p_1)-\cdots - p_n\ln(p_n)\leq\frac{\log_2(n)}{\log_2(\ee)}=\ln(n)\leq n \ln(n)p_n,
 \]
 where in the last inequality we used that \ $p_n=\max_{i\in\{1,\ldots,n\}} p_i$ \ implying \ $p_n\geq 1/n$.
\proofend

\section{$\cL_n$ \ and \ $\cP_n$ \ are not quasi arithmetic means}\label{App_non_quasi_arit_means}

 Given an interval \ $I\subset \RR$, \  and \ $n\geq 2$, \ $n\in\NN$, \ a map \ $M\colon I^n\to I$ \ is said to be bisymmetric
 if it fulfils the following equation
\begin{equation}\label{E:bisymmetry}
M(M(x_{11},\dots,x_{1n}),\dots,M(x_{n1},\dots,x_{nn}))=M(M(x_{11},\dots,x_{n1}),\dots,M(x_{1n},\dots,x_{nn}))
\end{equation}
for every $x_{ij}\in I,\ i,j=1,\dots,n$.

\begin{Thm}
If \ $n\geq 2$, \ $n\in\NN$, \ then \ $\cL_n$ \ is not a quasi arithmetic mean.
\end{Thm}

\noindent{\bf Proof.}
Let \ $n\geq 2$, \ $n\in\NN$, \ be fixed.
On the contrary, let us suppose that \ $\cL_n$ \  is a quasi arithmetic mean.
Then it should satisfy the following bisymmetry equation
\begin{equation}\label{E:bisymmetry_for_L_n}
\cL_n(\cL_n(x_{11},\dots,x_{1n}),\dots,\cL_n(x_{n1},\dots,x_{nn}))
   =\cL_n(\cL_n(x_{11},\dots,x_{n1}),\dots,\cL_n(x_{1n},\dots,x_{nn}))
\end{equation}
for all $x_{11},\dots,x_{1n},\dots,x_{n1},\dots,x_{nn}>1$, \ see, e.g.,  M\"unnich et al.\  \cite{MunMakMok}.

 Step 1.
\ We check that \eqref{E:bisymmetry_for_L_n} yields that the function \ $F\colon (1,\infty)\times(1,\infty)\to\RR$,
 \begin{align}\label{help_L_n_F}
 F(x,y):=\frac{\sqrt[n-1]{x}\ln(y)+\sqrt[n-1]{y}\ln(x)}{\ln(xy)},\qquad x,y\in(1,\infty),
 \end{align}
  should be bisymmetric as well.
 Here we will use the following extension of \ $\cL_n$:
\[
  \widetilde\cL_n(x_1,x_2,1,\dots,1):=\lim_{\substack{x_i\downarrow 1\\ { i\in\{3,\ldots,n\} }}}\cL_n(x_1,\dots,x_n)
                                     = \frac{\sqrt[n-1]{x_2}\ln(x_1) +\sqrt[n-1]{x_1}\ln(x_2)} {\ln(x_1x_2)},
                                     \qquad x_1,x_2>1.
\]
Let \ $x_{ij}>1,\ i,j=1,\dots,n$.
\ Taking the  iterated limits $x_{ij}\to 1+,\ i,j\not\in\{1,2\}$
  (in an arbitrary order) of both sides of \eqref{E:bisymmetry_for_L_n}, we have
\begin{align*}
&\widetilde\cL_n\left(\frac{\sqrt[n-1]{x_{12}} \ln(x_{11})+\sqrt[n-1]{x_{11}} \ln(x_{12})}{\ln(x_{11}x_{12})},
              \frac{\sqrt[n-1]{x_{22}}{ \ln(x_{21})} +\sqrt[n-1]{x_{21}}{ \ln(x_{22}) }}{\ln(x_{21}x_{22})},1,\dots,1\right)\\
&\qquad =\widetilde\cL_n\left(\frac{\sqrt[n-1]{x_{21}}\ln(x_{11}) +\sqrt[n-1]{x_{11}}\ln(x_{21}) }{\ln(x_{11}x_{21})},
     \frac{\sqrt[n-1]{x_{22}}{ \ln(x_{12})} + \sqrt[n-1]{x_{12}}{ \ln(x_{22})} }{\ln(x_{12}x_{22})},1,\dots,1\right),
\end{align*}
 where we used that
 \begin{align*}
  &\lim_{x_1\downarrow 1} \cL_n(x_1,\ldots,x_n)
     = \frac{\sum_{i=2}^n \sqrt[n-1]{\prod_{j=2,j\ne i}^n x_j}\ln(x_i)  }{\sum_{i=2}^n \ln(x_i)},
     \qquad x_2,\ldots,x_n>1,\\
  &\phantom{\lim_{x_1\downarrow 1} \cL_n(x_1,\ldots,x_n)} \vdots\\
  &\lim_{x_{n-1}\downarrow 1}\cdots\lim_{x_1\downarrow 1} \cL_n(x_1,\ldots,x_n)
     = \sqrt[n-1]{x_n},\qquad x_n>1,\\
  & \lim_{x_n\downarrow 1}  \lim_{x_{n-1}\downarrow 1}\cdots\lim_{x_1\downarrow 1} \cL_n(x_1,\ldots,x_n)  = 1.
 \end{align*}
Introducing the notations
\[
x_{11} =:x,\qquad x_{12}=:y,\qquad x_{21} =:s,\qquad x_{22} =:t,
\]
and, using the  definitions of $\cL_n$ and $\widetilde\cL_n$, we get
 \begin{align}\label{E:bisymmetry_for_F}
 \begin{split}
  &\frac{\sqrt[n-1]{F(s,t)}\ln(F(x,y)) + \sqrt[n-1]{F(x,y)} \ln(F(s,t))}{\ln(F(x,y) F(s,t))}\\
  &\qquad = \frac{\sqrt[n-1]{F(y,t)}\ln(F(x,s)) + \sqrt[n-1]{F(x,s)}\ln(F(y,t))}{\ln(F(x,s) F(y,t))},
 \end{split}
\end{align}
 i.e.,
 \[
  F(F(x,y), F(s,t)) = F(F(x,s), F(y,t)), \qquad x,y,s,t>1,
 \]
 yielding that \ $F$ \ is bisymmetric.

 Step 2.
\ We check that the function \ $F$ \ defined in \eqref{help_L_n_F} is not bisymmetric.
{On the contrary, let} us assume that \ $F$ \ is bisymmetric,  i.e., \eqref{E:bisymmetry_for_F} holds for all \ $x,y,s,t>1$.
\ We  distinguish two cases, \ $n>2$ \ and \ $n=2$.

\noindent At first, let $n>2$.
 By substituting \ $x=y=\ee^{2(n-1)^2}$ and $s=t=\ee^{(n-1)^2}$ in \eqref{E:bisymmetry_for_F},  after some simplifications and rearrangements,
 we get that
 \[
 \frac{\ee}{3}(\ee+2)=\sqrt[n-1]{\frac{\ee^{n-1}}{3}(\ee^{n-1}+2)}.
 \]
 Since the function \ $(0,\infty)\ni z\mapsto z^{n-1}$ \ is strictly convex for all \ $n>2$, \ we have
 \[
 \left(\frac{\ee+2}{3}\right)^{n-1}<\frac{\ee^{n-1}+2}{3},
 \]
which entails that \ $F$ \ can not be bisymmetric  for \ $n>2$.

For the case \ $n=2$, \ let  us substitute \ $x=y,\ s = \ee$, and $t = \ee^2$ in \eqref{E:bisymmetry_for_F}.
Then we  get
\[
 \frac{\frac{\ee}{3}(\ee+2)\ln(x)+x\ln\left(\frac{\ee}{3}(\ee+2)\right)}{\ln\left(x\frac{\ee}{3}(\ee+2)\right)}
  =\frac{\frac{2x+\ee^2\ln(x)}{\ln(x)+2}\ln\left(\frac{x+\ee\ln(x)}{\ln(x)+1}\right)
       +\frac{x+\ee\ln(x)}{\ln(x)+1}\ln\left(\frac{2x+\ee^2\ln(x)}{\ln(x)+2}\right)}
        {\ln\left(\frac{2x+\ee^2\ln(x)}{\ln(x)+2}\cdot\frac{x+\ee\ln(x)}{\ln(x)+1}\right)}.
\]
If we calculate  the values of both sides of the previous equation with \ $x=\ee^{10}$, \
 then we get  strictly less than \ $2800$ \ for the left hand side  (approximately \ $2797.9$),
 \ and strictly greater than \ $2800$ \ for the right hand side  (approximately \ $2808.8$).
\ So \ $F$ \ can not be bisymmetric even for \ $n=2$.

 Steps 1 and 2  lead us to a contradiction.
\proofend

\begin{Thm}
If \ $n\geq 2$, \ $n\in\NN$, \ then \ $\cP_n$ \ is not a quasi arithmetic mean.
\end{Thm}

\noindent{\bf Proof.}
Let \ $n\geq 2$, \ $n\in\NN$, \ be fixed.
We divide the proof into three steps.

 Step 1.
\ We check that \ $\cP_n$ \ is a quasi arithmetic mean on \ $(1,\infty)$ \ if and only if \ $\widetilde{\cP}_n$ \ is
 a quasi arithmetic mean on \ $(0,\infty)$, \ where
 \begin{align}\label{help_P_n_1}
 \widetilde{\cP}_n(y_1,\dots,y_n):=\frac{1}{n\ln(n)}\sum\limits_{i=1}^{n}y_i\ln\left(\frac{y_1+\cdots+y_n}{y_i}\right),\qquad y_1,\dots,y_n> 0.
 \end{align}

\noindent First, let us assume that \ $\cP_n$ \ is a quasi arithmetic mean on \ $(1,\infty)$.
\ Then there exists a strictly monotone increasing, continuous function \ $\varphi\colon (1,\infty)\to\RR$ \ such that
\[
  \cP_n(x_1,\ldots,x_n) =
 \left( \prod_{i=1}^n x_i^{\ln\left( \frac{\ln(x_1\cdots x_n)}{\ln(x_i)} \right)}\right)^{\frac{1}{n\ln(n)}}=\varphi^{-1}\left(\frac{\varphi(x_1)+\cdots+\varphi(x_n)}{n}\right),\qquad x_1,\dots,x_n> 1.
\]
With the substitutions
\begin{align}\label{help_P_n_2}
\ln(x_i) =: y_i,\quad i=1,\dots,n,\qquad\varphi\circ\exp =:f,
\end{align}
we can derive the equation
 \[
 \widetilde{\cP}_n(y_1,\dots,y_n)=f^{-1}\left(\frac{f(y_1)+\cdots+f(y_n)}{n}\right),\qquad y_1,\dots,y_n>0,
 \]
 yielding that \ $\widetilde{\cP}_n$ \ is a quasi arithmetic mean on \ $(0,\infty)$ \ corresponding to \ $f$.

\noindent Let us assume now that \ $\widetilde{\cP}_n$ \ is a quasi arithmetic mean on \ $(0,\infty)$.
\ Then there exists a strictly monotone increasing, continuous function \ $f:(0,\infty)\to\RR$ \ such that
 \[
   \widetilde{\cP}_n(y_1,\dots,y_n)
     = \frac{1}{n\ln(n)}\sum\limits_{i=1}^{n}y_i\ln\left(\frac{y_1+\cdots+y_n}{y_i}\right)
     = f^{-1}\left(\frac{f(y_1)+\cdots+f(y_n)}{n}\right)
 \]
 for \ $y_1,\ldots,y_n>0$.
\ With the substitutions \eqref{help_P_n_2}, we have
 \[
   \cP_n(x_1,\ldots,x_n) = \varphi^{-1}\left(\frac{\varphi(x_1)+\cdots+\varphi(x_n)}{n}\right),
      \qquad x_1,\ldots,x_n>1,
 \]
 yielding that \ $\cP_n$ \ is a quasi arithmetic mean on \ $(1,\infty)$ \ corresponding to \ $\varphi$.

Step 2.
\ We check that if $\widetilde{\cP}_n$ given in \eqref{help_P_n_1} is bisymmetric, then the
 function $ G:(0,\infty)\times (0,\infty)\to\RR$,
 \begin{align}\label{help_P_n_3}
   G(a,b):=\ln\left(\frac{(a+b)^{a+b}}{a^ab^b}\right),\qquad a,b>0,
 \end{align}
 should be bisymmetric as well.

If \ $\widetilde{\cP}_n$ \ is bisymmetric, then it fulfils the bisymmetry equation
\begin{equation}\label{E:bisymmetry_for_P_n}
\widetilde{\cP}_n(\widetilde{\cP}_n(y_{11},\dots,y_{1n}),\dots,\widetilde{\cP}_n(y_{n1},\dots,y_{nn}))
   =\widetilde{\cP}_n(\widetilde{\cP}_n(y_{11},\dots,y_{n1} ),\dots,\widetilde{\cP}_n(y_{1n},\dots,y_{nn}))
\end{equation}
for all $y_{11},\dots,y_{1n},\dots,y_{n1},\dots,y_{nn}>0$.
 Here we will use the following extension of $\widetilde{\cP}_n$:
\begin{align*}
  \widetilde{\cP}^*_n(y_1,y_2,0,\dots,0)
    &:=\lim_{\substack{y_i\downarrow 0\\  i\in\{3,\ldots,n\}}}\widetilde{\cP}_n(y_1,\dots,y_n)
      = \frac{1}{n\ln(n)} \left( y_1\ln\left(\frac{y_1+y_2}{y_1}\right)+y_2\ln\left(\frac{y_1+y_2}{y_2}\right)\right)\\
     &\;= \frac{1}{n\ln(n)} G(y_1,y_2),\qquad y_1,y_2>0.
\end{align*}
Let \ $ y_{ij}>0,\ i,j=1,\dots,n$.
\ Taking the  iterated limits $y_{ij}\downarrow 0,\ i,j\not\in\{1,2\}$
  (in an arbitrary order) of both sides of \eqref{E:bisymmetry_for_P_n}, we have%
 {\scriptsize
\begin{align*}
 &\widetilde{\cP}^*_n\left(\frac{1}{n\ln(n)}\left(y_{11}\ln\left(\frac{y_{11}+y_{12}}{y_{11}}\right)+y_{12}\ln\left(\frac{y_{11}+y_{12}}{y_{12}}\right)\right),
      \frac{1}{n\ln(n)}\left(y_{21}\ln\left(\frac{y_{21}+y_{22}}{y_{21}}\right)+y_{22}\ln\left(\frac{y_{21}+y_{22}}{y_{22}}\right)\right),0,\dots,0\right)\\
 &\quad =\widetilde{\cP}^*_n\left(\frac{1}{n\ln(n)}\left(y_{11}\ln\left(\frac{y_{11}+y_{21}}{y_{11}}\right)+y_{21}\ln\left(\frac{y_{11}+y_{21}}{y_{21}}\right)\right),
     \frac{1}{n\ln(n)}\left(y_{12}\ln\left(\frac{y_{12}+y_{22}}{y_{12}}\right)+y_{22}\ln\left(\frac{y_{12}+y_{22}}{y_{22}}\right)\right),0,\dots,0\right),
\end{align*}
}
 where we used that
 \begin{align*}
  &\lim_{y_1\downarrow 0} \widetilde{\cP}_n(y_1,\ldots,y_n)
     = \frac{1}{n\ln(n)} \sum_{i=2}^{n-1} y_i\ln\left(\frac{ y_2+\cdots+y_{n-1}}{y_i}\right),
     \qquad  y_2,\ldots,y_{n-1}>0,\\
  &\phantom{\lim_{x_1\downarrow 1} \cP_n(x_1,\ldots,x_n)} \vdots\\
  & \lim_{y_{n-2}\downarrow 0} \cdots\lim_{y_1\downarrow 0}   \widetilde{\cP}_n(y_1,\ldots,y_n)
      =\frac{1}{n\ln(n)} \!\left[y_{n-1}\ln\!\left(\!\frac{y_{n-1}+y_n}{y_{n-1}}\right)+y_n\ln\!\left(\!\frac{y_{n-1}+y_n}{y_n}\right)\right],\; y_{n-1},\,y_n>0,\\
  &  \lim_{y_{n-1}\downarrow 0} \cdots \lim_{y_1\downarrow 0}\widetilde{\cP}_n(y_1,\ldots,y_n)
       =\lim_{y_{n}\downarrow 0} \cdots \lim_{y_1\downarrow 0}\widetilde{\cP}_n(y_1,\ldots,y_n)=0,
       \qquad  y_n>0.
 \end{align*}
Introducing the notations
\[
y_{11} =:x,\qquad y_{12}=:y,\qquad y_{21} =:s,\qquad y_{22} =:t,
\]
and, using the definitions of $\widetilde{\cP}_n$ and $\widetilde{\cP}^*_n$, after some simplification, we get
 \begin{align*}
  & G(x,y)\ln\left(\frac{ G(x,y)+ G(s,t)}{ G(x,y)}\right)+ G(s,t)\ln\left(\frac{ G(x,y)+ G(s,t)}{ G(s,t)}\right)\\
  &\qquad  =G(x,s)\ln\left(\frac{G(x,s)+G(y,t)}{G(x,s)}\right)+G(y,t)\ln\left(\frac{G(x,s)+G(y,t)}{G(y,t)}\right),
\end{align*}
 i.e.,
 \[
  G(G(x,y), G(s,t)) = G(G(x,s), G(y,t)), \qquad x,y,s,t>0,
 \]
 yielding that \ $G$ \ is bisymmetric.

Step 3.
\ We check that the function \ $G$ \ defined in \eqref{help_P_n_3} is not bisymmetric, yielding that
 \ $\widetilde{\cP}_n$ \ can not be bisymmetric (see Step 2), and hence \ $\cP_n$ \ can not be bisymmetric (see Step 1).
On the contrary, let us suppose that \ $G$ \ is bisymmetric.
Note that \ $G$ \ is  strictly monotone increasing  in both of its variables, and continuous as well.
Hence according to Maksa \cite[Theorem 1]{Maksa1999}, there exist strictly monotone, continuous functions
 \  $\varphi_1\colon (0,\infty)\to\RR$, \ $\varphi_2\colon (0,\infty)\to\RR$, \ and
 \ $\psi\colon G((0,\infty)\times (0,\infty))\to\RR$ \ such that
 \[
 G(a,b)=\psi^{-1}( \varphi_1(a)+\varphi_2(b)), \qquad a,b\in(0,\infty).
 \]
 Since \ $G$ \ is symmetric as well, we have \ $\varphi_1(a)+\varphi_2(b) = \varphi_1(b)+\varphi_2(a)$, \ $a,b\in(0,\infty)$,
 \ i.e., \ $(\varphi_1 - \varphi_2)(a) = (\varphi_1 - \varphi_2)(b)$, \ $a,b\in(0,\infty)$, \ yielding the existence
 of \ $K\in\RR$ \ such that \ $\varphi_2(a) = \varphi_1(a) + K$, \ $a\in(0,\infty)$.
\ Hence
 \[
   G(a,b) = \psi^{-1}(\varphi_1(a) + \varphi_1(b) + K)
          = \psi^{-1}(\varphi(a) + \varphi(b)),
          \qquad a,b\in(0,\infty),
 \]
 where \ $\varphi:(0,\infty)\to\RR$, \ $\varphi(x):=\varphi_1(x)+\frac{K}{2}$, \ $x\in(0,\infty)$.
\ That is to say, \ $G$ \ is a quasisum  in the sense of Maksa \cite[Definition, page 59]{Maksa1999}.
Moreover, with the notation \ $h\colon (0,\infty)\to\RR,\ h(a):=a\ln (a)$, \ $a\in(0,\infty)$, \ $G$ \ can be written as a Cauchy-difference
 \[
 G(a,b)=h(a+b)-h(a)-h(b), \qquad a,b\in(0,\infty).
 \]
Since  \ $G$ \ is a Cauchy-difference and a quasisum at the same time, by J\'arai et al.\  \cite[Theorem 2.4]{JaraiMaksaPales2004},
 there exist an additive function \ $A\colon\mathbb{R}\to\mathbb{R}$ \ and \ $\alpha,\beta,\gamma,\delta\in\RR$ \ such that
 \ $\alpha\beta\not=0$ \ and \ $h$ \ should have one of the following forms:
\begin{enumerate}[(I)]
\item $h(x)=\alpha\ln(\cosh(\beta x+\gamma))+A(x)+\delta$;
\item $h(x)=\alpha\ln(\sinh(\beta x+\gamma))+A(x)+\delta$ \ (here  $\beta,\gamma\in\RR_+$);
\item $h(x)=\alpha\ln(\sin(\beta x+\gamma))+A(x)+\delta$ \  (here, in fact, \ $\beta=0$ \ and \ $\gamma\in(0,\pi)$);
\item $h(x)=\alpha \ee^{\beta x}+A(x)+\delta$;
\item $h(x)=\alpha\ln(|x+\gamma|)+A(x)+\delta$ \ (here  $\gamma\in\RR_+$);
\item $h(x)=\alpha x^2+A(x)+\delta$
\end{enumerate}
for all $x\in (0,\infty)$.
\ Since a continuous and additive function on \ $ (0,\infty)$ \ has the form \ $cx$, \ $x\in  (0,\infty)$,
 \ with some \ $c\in\RR$, \ it is clear that all the cases \ (I)--(VI) \ are impossible, so \ $G$ \ can not be bisymmetric.

Steps 1, 2 and 3 imply the assertion, since if \ $\cP_n$ \ were a quasi arithmetic mean,
 then it should be bisymmetric (see, e.g.,  M\"unnich et al.\  \cite{MunMakMok}),  leading us to a contradiction.
\proofend

\section{Application of means to congressional apportionment in the USA's election}\label{Ex_Toredek_szavazatok}

In the USA, the membership of the House of Representatives is fixed at 435  by the Apportionment Act of 1911,
 and the representation of each state in the House of Representatives is based on its population.
In principle, it would mean that the number of representatives of a given state in the House of Representatives
  can be calculated as follows: we multiply \ 435 \ by the population of the given state and divide it by the total population of the USA.
However, this number is not an integer in general, so, in practice, its integer part is taken (if it is \ $0$, \ then
 \ $1$ \ representative is apportioned to the given state).
As a result of this procedure there are some remaining places for representatives which should be apportioned among the  50 states.
This is an  important question, since there is a census in the USA in every $10^{\mathrm{th}}$ year (the next one will be in 2020).
Sullivan \cite{Sul1, Sul2} provided several methods for the apportionment such as the method of the arithmetic, geometric and harmonic means.
In what follows, we provide a common generalization of these three methods to quasi arithmetic means, and we also
 point out further possible extensions to Bajraktarevi\'c means and Cauchy quotient means.

 Let \ $N_A$ \ and \ $N_B$ \ be the population size of two states \ $A$ \ and \ $B$ \ in the USA, respectively,
 and \ $r_A$ \ and \ $r_B$ \ be the corresponding number of representatives assigned these states.
 Ideally, the ratios \ $\frac{r_A}{N_A}$ \ and \ $\frac{r_B}{N_B}$
 should be equal, however, in reality, this is not the case.
According to Sullivan's arithmetic method, one says that the assignment of an additional representative to  state \ $A$ \ rather than to  state \ $B$ \  is correct (fair) if
\begin{align*}
   \frac{r_A+1}{N_A} - \frac{r_B}{N_B}
      < \frac{r_B+1}{N_B} - \frac{r_A}{N_A},
 \end{align*}
  or equivalently
 \[
     \frac{1}{2}\left(\frac{r_A}{N_A} + \frac{r_A+1}{N_A}\right) < \frac{1}{2}\left(\frac{r_B}{N_B} + \frac{r_B+1}{N_B}\right),
 \]
 see Sullivan \cite{Sul1}.
Then one can arrange the values \ $\frac{1}{2} (\frac{r_i}{N_i} + \frac{r_i+1}{N_i})$, \ $i=1,\ldots,50$, \ in an increasing order,
 where \ $N_1,\ldots,N_{50}$ \ are the populations of the 50 states and \ $r_1,\ldots,r_{50}$ \ are the corresponding number of representatives
 (before assigning the remaining places).
If there are \ $k$ \ remaining places for representatives, then assign a representative  to those \ $k$ \ states which correspond to
 the bottom \ $k$ \ values in the above mentioned list.

 Let \ $f:(0,\infty)\to \RR$ \ be a continuous and strictly monotone increasing function.
The ratios  $f(r_A/N_A)$ \ and \ $f(r_B/N_B)$ \  are, as before, not equal in general.
Analogously to Sullivan's fairness definition, we say that the assignment of an additional representative
 to  state \ $A$ \ rather than to  state \ $B$ \  is fair with respect to the function $f$ if
 \begin{align*}
   f\left(\frac{r_A+1}{N_A}\right) - f\left(\frac{r_B}{N_B}\right)
      < f\left(\frac{r_B+1}{N_B}\right) - f\left(\frac{r_A}{N_A}\right),
 \end{align*}
 or equivalently
 \begin{align}\label{help_tor_szav1}
  M_2^f\left( \frac{r_A}{N_A}, \frac{r_A+1}{N_A} \right)
     < M_2^f\left( \frac{r_B}{N_B}, \frac{r_B+1}{N_B} \right),
 \end{align}
  where \ $M_2^f$ \ is the \ $2$-variable quasi arithmetic mean corresponding to \ $f$.
\ By choosing \ $f:(0,\infty)\to \RR$, \ $f(x)=x$, \ $f(x)=\ln(x)$ \ and \ $f(x)=x^{-1}$, \ $x>0$, \ one gets back the method of
 arithmetic, geometric and harmonic mean, respectively, given in Sullivan \cite{Sul1}.
Then one can arrange the values \ $M_2^f(\frac{r_i}{N_i}, \frac{r_i+1}{N_i})$, \ $i=1,\ldots,50$, \ in an increasing order, and, similarly as in the case of
 Sullivan's arithmetic method, if there are \ $k$ \ remaining places for representatives, then assign a representative  to those \ $k$ \ states which correspond to
 the bottom \ $k$ \ values in the above mentioned list.
As a generalization, one may replace the inequality \eqref{help_tor_szav1} by
 \begin{align}\label{help_tor_szav2}
  B_2^{f,p} \left( \frac{r_A}{N_A}, \frac{r_A+1}{N_A} \right)
     < B_2^{f,p}\left( \frac{r_B}{N_B}, \frac{r_B+1}{N_B} \right),
 \end{align}
 where \ $p:(0,\infty)\to\RR$ \ is a given (weight) function,
  where \ $B_2^{f,p}$ \ is the \ $2$-variable Bajraktarevi\'c mean corresponding to \ $f$ \ and \ $p$.
\ By some algebraic calculations, one can check that \eqref{help_tor_szav2} is equivalent to
 \begin{align*}
  &p\left(\frac{r_A+1}{N_A}\right)p\left(\frac{r_B+1}{N_B}\right)
       \left[f\left(\frac{r_B+1}{N_B}\right) - f\left(\frac{r_A+1}{N_A}\right) \right] \\
  &+ p\left(\frac{r_A+1}{N_A}\right)p\left(\frac{r_B}{N_B}\right)
       \left[f\left(\frac{r_B}{N_B}\right) - f\left(\frac{r_A+1}{N_A}\right) \right] \\
  &+ p\left(\frac{r_A}{N_A}\right)p\left(\frac{r_B+1}{N_B}\right)
       \left[f\left(\frac{r_B+1}{N_B}\right) - f\left(\frac{r_A}{N_A}\right) \right] \\
  &+ p\left(\frac{r_A}{N_A}\right)p\left(\frac{r_B}{N_B}\right)
       \left[f\left(\frac{r_B}{N_B}\right) - f\left(\frac{r_A}{N_A}\right) \right]
   > 0.
 \end{align*}
If one replaces the inequality \eqref{help_tor_szav1} by
 \begin{align}\label{help_tor_szav3}
    \cB_2\left( \frac{r_A}{N_A}, \frac{r_A+1}{N_A} \right)
       < \cB_2 \left( \frac{r_B}{N_B}, \frac{r_B+1}{N_B} \right),
 \end{align}
 then one gets back the method of harmonic mean in Sullivan \cite{Sul1},
 since \ $\cB_2(x,y)$ \ is nothing else but the harmonic mean of \ $x,y\in(0,\infty)$, and it is easy to check that
 \eqref{help_tor_szav3} is equivalent to
 \[
   \frac{N_A}{r_A+1} - \frac{N_B}{r_B} > \frac{N_B}{r_B+1} - \frac{N_A}{r_A}.
 \]
In general, in the inequality \eqref{help_tor_szav1} the quasi arithmetic mean \ $M_2^f$ \ corresponding to \ $f$ \
 could be replaced by any $2$-variable symmetric mean, and one could investigate the effects of the corresponding assignment rules
 for a given election in the USA, similarly as in Sullivan \cite{Sul1, Sul2}.

\section*{Acknowledgements}

The paper of Himmel and Matkowski \cite{HimMat3} was not publicly available when writing our paper,
 and we would like to thank Mih\'aly Bessenyei for asking Janusz Matkowski about the paper \cite{HimMat3}.
Janusz Matkowski was so kind  as to send Mih\'aly Bessenyei the paper in question, which finally reached us as well.
 We would like to thank the referee and the language editor for their comments that helped us to improve the paper.

\end{document}